\documentstyle{article}

\newcommand{\kl }{{\lambda}}

\newcommand{\ot }{\otimes}
\begin{document} 

\newcommand{\nt}{\noindent}

\title{Two Proofs of a Conjecture of Hori and Vafa}
\author{Aaron Bertram, Ionu\c t Ciocan-Fontanine and Bumsig Kim} 

\maketitle

\nt {\bf \S 0. Introduction.} The (small) quantum cohomology $QH^*(X)$ 
of a complex projective manifold $X$ can be thought of either as a 
deformation of the even degree cohomology
$V := H^{2*}(X,{\bf C})$ or as a family of
associative, commutative products on $V$ indexed by
$H^2(X,{\bf C})$. Indeed, if 
$T_1,...,T_m \in H^2(X,{\bf C})$ are a basis,
which we assume for simplicity to consist of nef divisors, 
(with dual basis $t_1,...,t_m$), then 
the small quantum product is a map:
$$*:V\times V \rightarrow V[[e^{t_1},...,e^{t_m}]]$$
(reducing to the cup product when $t_i \rightarrow -\infty$) determined by the
enumerative geometry of rational curves 
(the genus-zero Gromov-Witten invariants) on
$X$. If 
$X$ is a Fano manifold then $*$ is polynomial-valued in the exponentials. 

\medskip

The complex Grassmannian ${\bf G} := G(r,n)$ has one of the 
most-studied and best-understood
small quantum cohomology rings. Recall that the ordinary cohomology of 
${\bf G}$ 
has presentation:
$$H^*({\bf G}) \cong {\bf C}[\sigma_1,...,\sigma_r]/
\langle h_{n-r+1},...,h_n\rangle$$ 
where the $\sigma_i$ are the elementary symmetric polynomials in 
degree $2$ variables $x_1,...,x_r$
(the Chern roots of the dual to the universal bundle) 
and the $h_i$ are the complete
symmetric polynomials (sums of all monomials of 
degree $i$) in $x_1,...,x_r$. 

\medskip

The quantum cohomology of ${\bf G}$, on the other hand, has presentation:
$$QH^*({\bf G}) \cong {\bf C}[\sigma_1,...,\sigma_r,e^t]/
\langle h_{n-r+1},...,h_n + (-1)^re^t\rangle$$ 
where $T = \sigma_1$ is the basis for $H^2({\bf G},{\bf C})$ 
(see \cite{Wit}, \cite{ST}).

\medskip

Quantum cohomology has a mathematical partner that is 
frequently better suited for applications. Here one 
regards $*$ as an 
${\cal O}$-linear product 
of $TV$-valued vector fields over $H^2(X,{\bf C})$. That is, 
if $T_0,T_1,...,T_n$
is a basis for $V$, extending the basis of $H^2(X,{\bf C})$, 
with $T_0 = 1$ and $T_n$ its
Poincar\'e dual, and
if:
$$T_i * T_j = \sum_k \Phi^k_{ij}(e^{t_1},...,e^{t_m}) T_k$$
then 
$$*:TV|_{H^2} \otimes TV|_{H^2} \rightarrow TV|_{H^2}$$
is defined by $\partial_i * \partial_j = 
\sum_k \Phi^k_{ij}(e^{t_1},...,e^{t_m})  \partial_k$, 
where 
$\partial_i := \frac{\partial}{\partial t_i}$.


If we reinterpret $*$ once more as:
$$*:TV|_{H^2} \rightarrow TV|_{H^2} \otimes T^*V|_{H^2}$$
then the ${\cal O}$-linearity means that there is a family of connections:
$$\nabla_\hbar: TV|_{H^2} \rightarrow TV|_{H^2} 
\otimes T^*V|_{H^2}; \nabla_\hbar = 
d + \frac 1{\hbar} *$$
deforming $d$, and
the associativity of $*$ translates into the flatness of the $\nabla_\hbar$, 
i.e. existence of flat sections $F(t,\hbar) = F_0(t,\hbar)\partial_0 + ... +
F_n(t,\hbar)\partial_n$ satisfying:
$$\hbar \frac{\partial F}{\partial t_i} = 
\partial_i * F \ \mbox{for $i = 1,...,m$}$$
(see \cite{Dub}).

\medskip

Givental \cite{Giv1} computed the fundamental $(n+1)\times (n+1)$
matrix of solutions in terms of Gromov-Witten invariants 
``with gravitational descendents.'' The solution is unique of
the following form with suitable initial conditions:
$$(F_{i,j}(t,e^t,\hbar)\partial_j)_{0\leq i,j\leq n}$$
(polynomial in the $t = (t_1,...,t_m)$).
It is very  useful to regard the columns as
cohomology-valued, and the last column in particular:
$$J^X = \sum_{i=0}^n F_{i,n}(t,e^t,\hbar)\check{T}_i$$
is very functorial and plays an important role
in mirror symmetry (here $\{\check{T}_i\}$ is the basis dual
to $\{T_i\}$ with respect to the intersection pairing on $X$).

\medskip

In the simplest case, let $x$ be the hyperplane class in 
$H^2({\bf P}^{n-1})$. Then:
$$J^{{\bf P}^{n-1}} = e^{\frac{tx}{\hbar}}\sum_{d=0}^\infty \frac {e^{dt}}
{\prod_{l=1}^d (x +
l\hbar)^n} \ (\mbox{mod} \ x^n).$$
Our goal here is to compute 
the $J$-function of the Grassmannian ${\bf G}$ and explore various 
applications that result.
One of  the nice functorial properties of $J$-functions guarantees that:
$$J^{\bf P} = \prod_{i=1}^r J^{{\bf P}^{n-1}} = 
e^{\frac{t_1x_1 + ... + t_rx_r}{\hbar}} \sum_{(d_1,...,d_r)} 
\frac{e^{d_1t_1 + ... + d_rt_r}}
{\prod_{i=1}^r\prod_{l=1}^{d_i}(x_i + l\hbar)^n}$$ 
where ${\bf P} = \prod_{i=1}^r {\bf P}^{n-1}$ with hyperplane 
classes $x_i$ on the factors.

We will give two proofs of the following additional property of $J$:

\medskip

\nt {\bf Conjecture (Hori-Vafa,} [HV, Appendix A]{\bf ):} Let:
$$\Delta = \prod_{i< j} (x_i - x_j) \ \ \mbox{and}\ \ {\cal D}_\Delta = 
\prod_{i < j} 
\left(\hbar \frac{\partial}{\partial t_i} - 
\hbar \frac{\partial}{\partial t_j}\right)$$ 
denote the Vandermonde determinant and ``Vandermonde operator'', respectively. 
Then:
$$J^{\bf G} = \left(e^{-\sigma_1(r-1)\pi \sqrt{-1}/\hbar}\right)
\frac{{\cal D}_\Delta (J^{\bf
P})}{\Delta}|_{t_i = t + (r-1)\pi \sqrt{-1}}$$
i.e. the $J$-function of the Grassmannian is obtained by applying 
${\cal D}_\Delta$ 
to $J^{\bf P}$ and
then ``symmetrizing'' (and translating). 

\medskip

We will see in the
second proof of the conjecture, in particular, that 
this can be thought of as another nice
functorial property of
$J$-functions, which looks as though it ought to generalize 
to wider classes of geometric
invariant theory quotients.

We should note that Hori and Vafa formulate the conjecture in terms of
``period'' integrals, rather than $J$-functions. The connection between
the coefficients of the $J$-function and these integrals is explained
in Section 3. As an application, we then use the (equivariant) integral 
representation form of the Hori-Vafa conjecture to give a proof of Givental's
``$R$-Conjecture'', and hence of the Virasoro conjecture for Grassmannians.
Finally, we give (with Dennis Stanton) a proof for $G(2,n)$ of another
formula for one of the coefficients of the $J$-function, which
was conjectured in \cite{BCKS1} from different considerations.

\medskip

\nt {\bf Acknowledgements:} The authors are aware of at least two 
previous attempts to compute the $J$-function of the Grassmannian 
(or an equivalent formulation). In particular,
recent work of \cite{LLLY} contains Lemmas 1.1 and 1.2 of 
this paper, and then produces
a lengthy algorithm, but not a closed formula, for the 
coefficients of the $J$-function.
Our first proof of  the conjecture avoids many of the technical 
difficulties of that paper (and
of an earlier paper of \cite{Kim1} along similar lines) by using Euler 
sequences rather than relying 
on the equivariant cohomology ring of the Grassmannian. 


Our second proof relies upon a quantum cohomology
interpretation of the Hori-Vafa conjecture that was suggested to
the second author by Sheldon Katz in the Fall of 2000. 


We thank Dennis Stanton for providing the proof of Proposition 3.5,
and for allowing us to include it in this paper. We also thank
Alexander Givental, Dosang Joe, J. Park, and Alexander Yong for
useful discussions.

A. Bertram and I. Ciocan-Fontanine have been partially supported
by NSF grants DMS-0200895 and DMS-0196209, respectively. 
B. Kim has been 
supported by KOSEF 1999-2-102-003-5, R03-2001-00001-0,
R02-2002-00-00134-0.

\medskip

\newpage

\nt {\bf \S 1. A Localization Proof of the Hori-Vafa Conjecture.} Let: 
$${\rm Quot}_{{\bf P}^1,d}({\bf C}^n,n-r)$$
be the Grothendieck quot scheme parametrizing the coherent-sheaf quotients 
${\bf C}^n\otimes{\cal O}_{{\bf P}^1} \rightarrow {\cal Q}$ with 
Hilbert polynomial $d + (n - r)(t+1)$ 
(see e.g. \cite{BDW}) (this is the Hilbert polynomial
of a locally free sheaf on ${\bf P}^1$ of rank $n- r$ and degree $d$). 
The quot scheme  
is a smooth, projective variety and if $n > r$, it is a 
compactification of the Hilbert
scheme of maps $g:{\bf P}^1 \rightarrow {\bf G}$ of degree $d$
since a {\it locally free} quotient of 
${\bf C}^n\otimes {\cal O}_{{\bf P}^1}$ of rank $n-r$ and degree $d$ is a
map to ${\bf G}$. 

\medskip

Consider on the other hand the quot scheme:
$${\rm Quot}_{{\bf P}^1,d}({\bf C}^r,0)$$
of {\it torsion} quotients of ${\bf C}^r\otimes {\cal O}_{{\bf
P}^1}$ of length
$d$. These were considered by Weil, as higher rank versions
of the symmetric product
${\bf P}^d = {\rm Sym}^d{\bf P}^1 = {\rm Quot}_{{\bf P}^1,d}({\bf C},0)$
Indeed, every such quot scheme maps to the symmetric product:
$$\wedge^r: {\rm Quot}_{{\bf P}^1,d}({\bf C}^r,0) \rightarrow
{\rm Quot}_{{\bf P}^1,d}({\bf C},0) = {\bf P}^d$$
via the top exterior power of the {\it kernel} vector bundle:
$$\wedge^r(K \subset {\bf C}^r\otimes {\cal O}_{{\bf P}^1}) := 
(\wedge^rK \subset {\bf C}\otimes {\cal O}_{{\bf P}^1})$$

Each such kernel is a locally free sheaf on ${\bf P}^1$ with a {\it splitting
type}:
$$K \cong {\cal O}_{{\bf P}^1}(-d_1) \oplus ... \oplus 
{\cal O}_{{\bf P}^1}(-d_r)$$
such that $d_1 + ... + d_r = d$. This is unique when we require 
$0 \le d_1 \le d_2 \le ... \le d_r$.
The action of PGL$(2,{\bf C})$ on ${\bf P}^1$ determines an 
action on the quot schemes
by pulling back kernels. For such actions,
$\wedge^r$ is equivariant and the splitting type of the 
kernel $K$ is evidently invariant.

\medskip

We will consider the diagonal action of ${\bf C}^*
\subset \mbox{PGL}(2,{\bf C})$ on ${\bf P}^1$:
$$\sigma(\zeta,(x:y)) = (\zeta x:y)$$
with fixed points at  $0 = (0:1)$ and $\infty = (1:0)$.
If ${\bf C}^n\otimes {\cal O}_{{\bf P}^1} \rightarrow {\cal Q}$
is a fixed point for the induced action of ${\bf C}^*$ on
${\rm Quot}_{{\bf P}^1,d}({\bf C}^r,0)$, then it is immediate 
that the (reduced) support
of ${\cal Q}$ is contained in $\{0,\infty\}$; the fixed points
of the action on ${\bf P}^1$. What is less immediate is the following:

\medskip

\nt {\bf Lemma 1.1:} For each splitting type
$\{d_i\}$ as above, let 
$m_1,m_2,...,m_k$ denote the jumping indices
(i.e. $0\leq d_1 = ... = d_{m_1} < d_{m_1+1} =
... = d_{m_2} < ... $). Then there is an embedding:
$$i_{\{d_i\}}: Fl(m_1,m_2,...,m_k,r) 
\hookrightarrow {\rm Quot}_{{\bf P}^1,d}({\bf C}^r,0)$$ 
with the property that each fixed point of the ${\bf C}^*$-action with
supp$({\cal Q}) = \{0\}$ and with kernel splitting type 
$\{d_i\}$ corresponds to a
point of  the flag variety.

\medskip

{\bf Proof:} Start with the universal flag on $Fl :=
Fl(m_1,m_2,...,m_k,r)$
$$0 \subseteq S_{m_1} \subset S_{m_2} \subset ... \subset S_{m_k} \subset
S_{m_{k+1}} = {\bf C}^r\otimes {\cal O}_{Fl}$$ let $\pi: {\bf P}^1 \times Fl
\rightarrow Fl$ be the projection and let
$z = 0 \times Fl
\subset  {\bf P}^1\times Fl$. Then we construct a modified flag of 
subsheaves of ${\bf
C}^r\otimes {\cal O}_{{\bf P}^1\times Fl}$ as follows. First, 
construct $S'_{m_i}$
for $i > 1$ via an elementary modification:
$$\begin{array}{ccccccccc}
0 & \rightarrow & \pi^*S_{m_1} & \rightarrow & S'_{m_i} &
\rightarrow & \pi^*S_{m_i}/S_{m_1}(-d_{m_2}z) &
\rightarrow & 0 \\
& & \| & & \downarrow & & \downarrow \\
0 & \rightarrow & \pi^*S_{m_1} & \rightarrow & \pi^*S_{m_i} & \rightarrow &
\pi^*S_{m_i}/S_{m_1} &
\rightarrow & 0 \end{array}$$
and note that this gives a flag: $\pi^*S_{m_1} \subset S'_{m_2} 
\subset ... \subset
S'_{m_{k+1}} \subset {\bf C}^r\otimes {\cal O}_{{\bf P}^1\times Fl}$.
Then, inductively, define $S^{(j)}_{m_i}$ for $i > j$ by:

$$\begin{array}{ccccccccc}
0 & \rightarrow & S^{(j-1)}_{m_j} & \rightarrow & S^{(j)}_{m_i} &
\rightarrow & S^{(j-1)}_{m_i}/S^{(j-1)}_{m_j}(-(d_{m_{j+1}}- d_{m_j})z) &
\rightarrow & 0 \\
& & \| & & \downarrow & & \downarrow \\
0 & \rightarrow & S^{(j-1)}_{m_j} & \rightarrow & S^{(j-1)}_{m_i} &
\rightarrow & S^{(j-1)}_{m_i}/S^{(j-1)}_{m_j} &
\rightarrow & 0\end{array}$$

This process yields, in the end, a flag of sheaf-inclusions of vector bundles:

$$0 \subset \pi^*S_{m_1} \subset S^{(1)}_{m_2} \subset ... \subset
S^{(i-1)}_{m_i} \subset ... \subset S^{(k)}_{m_{k+1}} \subset
{\bf C}^r\otimes {\cal O}_{{\bf P}^1\times Fl}$$
with the property that
$$S^{(i-1)}_{m_i}/S^{(i-2)}_{m_{i-2}} \cong
\pi^*\left(S_{m_i}/S_{m_{i-1}}\right)(-d_{m_i}z)$$
and we define ${\cal K} := S^{(k)}_{m_{k+1}}$ with $Fl$-valued quotient
$${\bf C}^r\otimes {\cal O}_{{\bf P}^1\times Fl} \rightarrow
{\cal Q} = {\bf C}^r\otimes {\cal O}_{{\bf P}^1\times Fl}/{\cal K}$$

This quotient is flat over $Fl$, of (relative) length $d$
supported on $z$.
If we restrict to a point $f$ of the flag variety, we get a quotient
${\bf C}^r \otimes {\cal O}_{{\bf P}^1} 
\rightarrow {\cal Q}|_{{\bf P}^1\times f}$
with kernel $K_f$ of splitting type
$\{d_i\}$.

\medskip

It follows from the global description of the kernel that 
$K_f \rightarrow {\bf
C}^r\otimes {\cal O}_{{\bf P}^1}$ can be
expressed in matrix (block) form as:
$$\left[\begin{array}{cccccccccc} A_1x^{d_{m_1}} & A_2x^{d_{m_2}} & \cdots &
A_{k+1}x^{d_r}
\end{array}\right]:K \rightarrow {\bf C}^r \otimes {\cal O}_{{\bf P}^1}$$
where each $A_i$ is an $r-m_{i-1}\times m_i-m_{i-1}$ (block) matrix of scalars,
augmented below by zeroes (we set $m_0 = 0$). The $A_i$ represent
the (un-modified) flag of subspaces corresponding to $f$, and the $x^{d_{m_i}}$
factors are produced by the elementary modifications.

\medskip

On the other hand, if ${\cal Q}$ is a quotient of length $d$ supported
at $0$ and a fixed point for the action of ${\bf C}^*$, it follows that
every entry in the matrix for the map 
$K \rightarrow {\bf C}^r\otimes {\cal O}_{{\bf
P}^1}$ is a multiple of a power of $x$. 
It follows that after an automorphism of
$K$, the matrix can be put into the block form above, hence is in the image of
$i_{\{d_i\}}$ where $\{d_i\}$ is the splitting type of
$K$. Finally, the ambiguity in the block form of the matrix is 
the same on both sides,
namely the action of $\times GL(m_i-m_{i-1},{\bf C})$. 
Thus the map $i_{\{d_i\}}$
is an embedding of smooth varieties and surjects onto the 
desired fixed loci. This
concludes the proof of  Lemma 1.1.

\bigskip

If $X$ is a scheme equipped with a vector bundle $E$ of rank $r$, 
then there is a
relative version of Lemma 1.1. Namely, the relative quot scheme over $X$:
$${\rm Quot}_{{\bf P}^1,d}(E,0) \rightarrow X$$
represents the functor: ``quotients
$\pi^*E \rightarrow {\cal Q}$
(for $\pi:{\bf P}^1\times T \rightarrow T\rightarrow X$) that are flat
of relative length $d$ over $T$''. By a theorem of Grothendieck, the 
fibers of the relative quot scheme over $X$ are isomorphic
to ${\rm Quot}_{{\bf P}^1,d}({\bf C}^r,0)$. Moreover, the ${\bf
C}^*$ action globalizes, and we obtain morphisms of $X$-schemes:
$$\begin{array}{ccc} i_{\{d_i\}}:Fl(m_1,...,m_k,E) & \longrightarrow & {\rm
Quot}_{{\bf P}^1,d}(E,0)
\end{array}$$ characterizing the fixed loci of the ${\bf C}^*$ action with
supp$({\cal Q}) = \{0 \times X\}$.

\medskip

We apply this version for the universal sub-bundle $S$ on $G(r,n)$ to get:

\medskip

\nt {\bf Lemma 1.2:} There is a natural ${\bf C}^*$-equivariant embedding:

$$i: {\rm Quot}_{{\bf P}^1,d}(S,0) \hookrightarrow 
{\rm Quot}_{{\bf P}^1,d}({\bf C}^n,n-r)$$
such that all the fixed points of the ${\bf C}^*$ action on
${\rm Quot}_{{\bf P}^1,d}({\bf C}^n,n-r)$ are contained in the image.
The fixed points of ${\rm Quot}_{{\bf P}^1,d}({\bf C}^n,n-r)$
that also satisfy supp$({\rm tor}({\cal Q})) = \{0\}$ 
(the support of the torsion part
of ${\cal Q}$) are precisely the images of flag manifolds:
$$Fl(m_1,m_2,...,m_k,r,n) = Fl(m_1,m_2,...,m_k,S)
\hookrightarrow {\rm Quot}_{{\bf P}^1,d}(S,0)$$
embedded by the relative version of Lemma 1.1.

\medskip

{\bf Proof:} The map $i$ is defined as follows. The kernel of the 
universal quotient:
$${\cal K} \hookrightarrow \pi^*S \rightarrow {\cal Q}$$
on ${\bf P}^1 \times {\rm Quot}_{{\bf P}^1,d}(S,0)$ can be 
thought of as a subsheaf of
${\bf C}^n\otimes {\cal O}_{{\bf P}^1 \times {\rm Quot}_{{\bf P}^1,d}(S,0)}$
by composing with the inclusion $\pi^*S \rightarrow 
\pi^*{\bf C}^n\otimes {\cal O}_G$.
Then the quotient:
$${\bf C}^n\otimes {\cal O}_{{\bf P}^1 \times {\rm Quot}_{{\bf
P}^1,d}(S,0)} \rightarrow {\cal Q}' = 
{\bf C}^n\otimes {\cal O}_{{\bf P}^1 \times {\rm
Quot}_{{\bf P}^1,d}(S,0)}/{\cal K}$$
is flat of the desired Hilbert polynomial. This gives the map $i$.

\medskip

The image of $i$ is the set of quotients ${\bf C}^n\otimes {\cal
O}_{{\bf P}^1} \rightarrow {\cal Q}'$ such that:
$$K \hookrightarrow {\bf C}^r\otimes {\cal O}_{{\bf P}^1} \rightarrow
{\bf C}^n\otimes {\cal O}_{{\bf P}^1}$$
with ${\cal Q} = {\bf C}^r\otimes {\cal O}_{{\bf P}^1}/K$ and ${\cal Q}' = 
{\bf C}^n\otimes {\cal O}_{{\bf P}^1}/K$.

\medskip

Now suppose that a quotient ${\bf C}^n \otimes {\cal O}_{{\bf P}^1} 
\rightarrow {\cal
Q}'$ in ${\rm Quot}_{{\bf P}^1,d}({\bf C}^n,n-r)$ is fixed under 
the action of ${\bf
C}^*$. Then the matrix of the map $K \rightarrow {\bf C}^n\otimes
{\cal O}_{{\bf P}^1}$ has $i$th column
$\vec v x^{b_i}y^{c_i}$ where $\vec v$ is a vector of scalars, and 
$b_i + c_i = d_i$.
From this it follows that $K$ factors through the subspace ${\bf C}^r 
\hookrightarrow
{\bf C}^n$ with matrix given by the $\vec v_i$. Finally, it is clear from the
construction that supp$({\rm {tor}}({\cal Q}')) = \mbox{supp}({\cal Q})$, hence
that the fixed loci with supp$({\rm {tor}(\cal Q}')) = \{0\}$ are the flag
manifolds of the global version on Lemma 1.1. This completes the 
proof of Lemma 1.2.

\medskip

Next, we need to relate these quot schemes to Kontsevich-Manin stacks. Recall
the:

\medskip

\nt {\bf Definition:} A map $f:C \rightarrow X$
from a curve $C$ with marked points $p_1,...,p_n \in C$ is {\it pre-stable} if:

\medskip

(i) $C$ is connected and projective, with at worst ordinary nodes 

\medskip

(ii) the marked points are smooth points of $C$

\medskip

\nt and $f$ is {\it stable} if, in addition:

\medskip

(iii) the automorphisms of $C$ fixing the $p_i$ and commuting with $f$ are
finite.
 
\medskip

Assuming that $X$ is homogeneous, there is a smooth proper stack
$\overline M_{0,n}(X,\beta)$
for each class $\beta \in H_2(X,{\bf Z})$ and each $n$ representing the functor
``flat families of stable maps of genus zero, 
$n$-pointed curves of class $\beta$.''  
Of special interest is the ``graph space:''
$$\overline M_{0,0}(X\times {\bf P}^1,(\beta,1))$$
which can be thought of as a compactification of the 
Hilbert scheme of maps $g:{\bf
P}^1\rightarrow X$ of class $\beta$ since a stable map 
$f:C \rightarrow X \times
{\bf P}^1$ of bidegree $(\beta,1)$ is the graph of a map
$g: {\bf P}^1 \rightarrow X$ whenever $C$ is irreducible. 
We will compare this compactification (when $X = {\bf G}$ and 
$\beta = d$) with 
the quot scheme. But first: 

\medskip

The ${\bf C}^*$ action on ${\bf P}^1$ induces an action on the graph spaces,
and one of the connected components of the fixed point locus is
$$\overline M_{0,1}(X,\beta) \cong F \subset 
\overline M_{0,0}(X\times {\bf P}^1,(\beta,1))$$
corresponding to the stable maps $f: C \rightarrow X 
\times {\bf P}^1$ with the property
that $C = C_0 \cup C_1, n = C_0 \cap C_1$ with $f(n) = (p,0)$,
$$f|_{C_0}:C_0 \rightarrow X 
\times \{0\} \ \mbox{of class $\beta$} \ \mbox{and}\ 
f|_{C_1}:C_1 \stackrel \sim \rightarrow p \times {\bf P}^1$$

This defines a stable map to $X\times {\bf P}^1$ exactly when the map $f_0$ 
is stable as a map from the pointed curve $(C_0,n)$. 

\medskip

\nt {\bf Theorem (Givental):} The  coefficients of the $J$-function of $X$
satisfy: 
$$J^X = e^{\frac {\sum_{i=1}^m t_iT_i}{\hbar}} \sum_{\beta} e^{\sum_{i=1}^m
t_i\int_{\beta}T_i } 
J_\beta(\hbar), $$
where $$J_\beta(\hbar) = ev_* \frac{1}{e_{{\bf
C}^*}(F)},$$ 
$e_{{\bf C}^*}(F)$ 
is the equivariant
Euler class of
$F$, computed in the ring:
$$H^*_{\bf C^*}(F,{\bf Q}) = H^*(F,{\bf Q}) \otimes _{\bf Q} {\bf Q}[\hbar]$$  
(i.e. we interpret $\hbar$ geometrically via $H^*(B{\bf C}^*,{\bf Q}) = 
{\bf Q}[\hbar]$) and 
$$ev:F = \overline M_{0,1}(X,\beta) \rightarrow X$$
is the evaluation map discussed above.

\medskip

Our task now is to compute this with the quot scheme when $X = {\bf G}$. 
For this, the following diagram (denoted $(\dag)$) is the key:

$$\begin{array}{ccccccc} {\rm Quot}_{{\bf P}^1,d}({\bf C}^n,n-r) &
\stackrel {\wedge^r} \rightarrow & {\bf P}^{\left(n \atop r\right)-1}_d &
\stackrel \Phi\leftarrow &
\overline M_{0,0}(G(r,n) \times {\bf P}^1,(d,1)) \\ \\
\cup & & \cup & &  \cup \\ \\
\coprod_{\{d_i\}} i_{\{d_i\}}(Fl) & \stackrel p\rightarrow & {\bf P}^{\left(n
\atop r\right)-1} &
\stackrel q\leftarrow &  F
\\ \\ & \searrow \rho & \cup & \ \ ev\swarrow \\ \\ & & G(r,n)
\end{array}$$

We need to explain this diagram. First, $G(r,n)$ is embedded by Pl\"ucker:
$$\wedge^r: G(r,n) \subset {\bf P}^{\left(n \atop r\right)-1}$$
Next, ${\bf P}^{\left(n \atop r\right)-1}_d = 
{\rm Quot}_{{\bf P}^1,d}({\bf C}^{\left(n
\atop r \right)},\left(n \atop r\right) - 1)$ is a projective space, 
more simply:
$${\bf P}^{\left(n \atop r\right)-1}_d = {\bf P}({\rm Hom}_d({\bf C}^2,{\bf
C}^{\left(n \atop r\right)}))$$
is the projectivized space of $d$-linear maps
${\rm Hom}_d({\bf C}^2,{\bf C}^n) = \mbox{Sym}^d({\bf C}^2)^*
\otimes {\bf C}^n$. 

\medskip

Thus the second row of $(\dag)$ consists of components of the fixed 
locus for the ${\bf
C}^*$-actions on the first row, with respect to which $\wedge^r$ 
and $\Phi$ are
equivariant. By Lemma 1.2, the morphisms $p$ and $q$ have, as their 
domains, all the fixed
points contained in
$(\wedge^r)^{-1}({\bf P}^{\left(n \atop r\right)-1})$ and 
$\Phi^{-1}({\bf P}^{\left(n \atop
r\right)-1})$ respectively. The restriction of $\rho$ to each 
component $i_{\{d_i\}}(Fl)$ is
the natural projection from the flag bundle 
$Fl(m_1,...,m_k,S)\rightarrow G(r,n)$.
It then follows
from the localization theorem of Atiyah-Bott (see \cite{Ber3}) that:

$$(\dag\dag)\ \ \wedge^r_*J_d(\hbar) = q_*\frac{1}{e_{{\bf C}^*}(F_0)} =
\sum_{\{d_i\}} p_*\left(\frac 1{e_{{\bf C}^*}(i_{\{d_i\}}(Fl))}\right)$$

We claim that more is true:

\medskip

\nt {\bf Lemma 1.3:} The $J$-function on $G(r,n)$ satisfies:
$$J_d(\hbar) = \sum_{\{d_i\}} 
\rho_*\left(\frac 1{e_{{\bf C}^*}(i_{\{d_i\}}(Fl))}\right)$$

{\bf Proof:} If the push-forward $\wedge^r_*:H^*(G(r,n),{\bf Q}) \rightarrow
H^*({\bf P}^{\left(n \atop r\right)-1},{\bf Q})$ were injective, 
this would follow
from $(*)$. Of course this is not the case, but the push-forward of
{\it equivariant cohomology rings} is injective. That is, the ``big''
torus $T = {{\bf C}^*}^n$ acting diagonally on ${\bf C}^n$, 
with $H^*(BT,{\bf Q})
= {\bf Q}[\lambda_1,...,\lambda_n]$ induces an action of the product
${\bf C}^*\times T$ on each of the spaces in $(\dag)$. All the
maps are equivariant for this product action, and it follows from  
the localization
theorem that the Pl\"ucker embedding induces an injective map:
$$\wedge^r_T:H_T^*(G(r,n),{\bf Q}) \rightarrow
H_T^*({\bf P}^{\left(n \atop r\right)-1},{\bf Q})$$
in the equivariant cohomology rings (for the $T$-action), and that moreover:
$$(\dag\dag)_T\ \  {q_T}_*\frac{1}{e_{{\bf C}^*\times T}(F_0)} =
\sum_{\{d_i\}} {p_T}_*\frac 1{e_{{\bf C}^*\times T}(i_{\{d_i\}}(Fl))} \in
H_T^*(G(r,n),{\bf Q})$$
for the ${\bf C}^*$-equivariant Euler classes with values in $T$-equivariant
cohomology, and $T$-equivariant push-forwards ${p_T}_*$ and ${q_T}_*$. Thus:

$${ev_T}_*\frac 1{e_{{\bf C}^*\times T}(F_0)} = 
\sum_{\{d_i\}} {\rho_T}_*\frac 1{e_{{\bf
C}^*\times T}(i_{\{d_i\}}(Fl))}$$
and since
$$\lim_{\lambda_i \rightarrow 0} {ev_T}_*\frac 1{e_{{\bf C}^*\times T}(F_0)} =
{ev}_*\frac 1{e_{{\bf C}^*}(F_0)} = J_d(\hbar)$$
and likewise for each 
${\rho_T}_*(1/{e_{{\bf C}^*\times T}(i_{\{d_i\}}(Fl))})$, 
the Lemma follows.

\bigskip

Thus we need to compute the equivariant Euler classes:
$$e_{{\bf C}^*}(i_{\{d_i\}}(Fl))$$
for each of the embeddings $i_{\{d_i\}}: Fl(m_1,...,m_k,S)
\rightarrow {\rm Quot}_{{\bf P}^1,d} ({\bf C}^n,n-r)$.

\medskip

Let:
$$x_{m_{i-1}+s} \ \ \mbox{for}\ \ i = 1,...,{k+1}\ \mbox{and}\ 
s = 1,...,m_{i}-m_{i-1}$$
denote the Chern roots of the universal vector bundles
$\left(S_{m_{i}}/S_{m_{i-1}}\right)^*$ on the flag bundle
(recall that in our notation $m_0=0$,
$m_{k+1}=r$, and $S_{m_{k+1}}=\rho ^*(S)$, with
$\rho :Fl(m_1,...,m_k,S)\rightarrow G(r,n)$ the projection).

\medskip

That is, the $x_{m_{i-1}+s}$ (formally) factor the Chern polynomials:
$$\prod_{s=1}^{m_i - m_{i-1}} (t + x_{m_{i-1}+s}) =
c_t(\left(S_{m_{i}}/S_{m_{i-1}}\right)^*) .$$
and in particular the $x_1,...,x_r$ are the Chern roots of $\rho ^*(S^*)$.

\medskip

Recall that
$S^*\otimes Q$ is the tangent bundle to the Grassmannian, where:
$$0 \rightarrow S \rightarrow {\bf C}^n\otimes {\cal O}_G \rightarrow Q 
\rightarrow 0$$
is the universal sequence of bundles on $G(r,n)$ and similarly,
$\pi_*({\cal K}^*\otimes {\cal Q})$ is the tangent bundle to
${\rm Quot}_{{\bf P}^1,d}({\bf C}^n,n-r)$, where
$$0 \rightarrow {\cal K} \rightarrow {\bf C}^n
\otimes {\cal O}_{{\bf P}^1\times {\rm
Quot}} \rightarrow {\cal Q} \rightarrow 0 $$
is the universal sequence of sheaves on ${\bf P}^1\times {\rm Quot}_{{\bf
P}^1,d}({\bf C}^n,n-r)$.

\medskip

From these, we obtain ``Euler sequences'':

$$0 \rightarrow S^*\otimes S \rightarrow S^*\otimes {\bf C}^n \rightarrow
TG(r,n) \rightarrow 0$$
and
$$0 \rightarrow \pi_*({\cal K}^*\otimes {\cal K}) \rightarrow 
\pi_*{\cal K}^*\otimes
{\bf C}^n \rightarrow T{\rm Quot} \rightarrow R^1\pi_*({\cal K}^*
\otimes {\cal K})
\rightarrow 0$$  

As for the flag bundle
$Fl(m_1,...,m_k,S)$, its tangent bundle
is given by:

$$0 \rightarrow K \rightarrow \rho^*S^*\otimes {\bf C}^n \rightarrow
TFl(m_1,...,m_k,S) \rightarrow 0$$
where $K$ has an increasing filtration $0 = K_0 \subset ... \subset K_k 
\subset K_{k+1} = K$  with
$K_{i}/K_{i-1} \cong (S_{m_{i}}/S_{m_{i-1}})^* \otimes S_{m_{i}}$
and the quotients $K_{i}/K_{i-1}$ filter further:
$$(S_{m_{i}}/S_{m_{i-1}})^* \otimes S_{m_1} \subset ... \subset
(S_{m_{i}}/S_{m_{i-1}})^*
\otimes S_{m_{i-1}} \subset (S_{m_{i}}/S_{m_{i-1}})^*
\otimes S_{m_{i}}$$
so that in the Grothendieck group of sheaves on $Fl(m_1,...,m_k,S)$:
$$[TFl] \sim
n[\rho^*S^*] - \sum_{i\le j}[(S_{m_{j}}/S_{m_{j-1}})^* \otimes
(S_{m_{i}}/S_{m_{i-1}})]$$

The restriction of ${\cal K}$ from ${\bf P}^1 \times \mbox{Quot}$ to 
${\bf P}^1\times Fl$ (also denoted by ${\cal K}$) also filters by
${\cal K}_i$ ($S^{(i-1)}_{m_i}$ in the proof of Lemma 1.1) with
$${\cal K}_{i}/{\cal K}_{i-1} \cong \pi^*(S_{m_{i}}/S_{m_{i-1}})(-d_{m_{i}}z)$$
so that in the Grothendieck group of sheaves on ${\bf P}^1\times Fl$:

$$[{\cal K}^*\otimes {\cal K}] \sim \sum_{i,j=1}^{k+1}
[\pi^*\left((S_{m_{i}}/S_{m_{i-1}})^*
\otimes (S_{m_{j}}/S_{m_{j-1}})\right)(d_{m_{i}} -
d_{m_{j}})]$$
and now we can compute 
$e_{{\bf C}^*}(i_{\{d_i\}}Fl) = 
c^{{\bf C}^*}_{\rm top}(T{\rm Quot}|_{i_{\{d_i\}}Fl}/ TFl)$
using
$$(**)\ \ \ c^{{\bf C}^*}_{\rm top}(T{\rm Quot}|_{i_{\{d_i\}}Fl}/TFl) = 
\frac{c^{{\bf C}^*}_{\rm
top}(\pi_*{\cal K}^*/\pi^*S^*)^n \cdot
c^{{\bf C}^*}_{\rm top}(R^1\pi_*({\cal K}^*\otimes {\cal K}))}
{c^{{\bf C}^*}_{\rm top}(\pi_*({\cal K}^*\otimes
{\cal K})/K)}$$ and replacing ${\cal K}$
and $K$ by their equivalent forms in the Grothendieck group. 

\medskip

Now if $E$ is a vector bundle of rank $e$ with 
with trivial action of ${\bf
C}^*$ and Chern roots $x_i$, then:

\medskip

(a) $c^{{\bf C}^*}_{\rm top}(H^0({\bf P}^1,{\cal O}_{{\bf P}^1}(d))\otimes E)=
\prod_{i=1}^e\prod_{l=0}^d(x_i + l\hbar)$

\medskip

\nt and if $y_j$ are the Chern roots of a bundle $F$ of rank $f$, then:

\medskip

(b) $c^{{\bf C}^*}_{\rm top}(H^0({\bf P}^1,{\cal O}_{{\bf P}^1}(d))
\otimes E^*\otimes F) =
\prod_{i=1}^e\prod_{j=1}^f\prod_{l=0}^d(y_j - x_i + l\hbar)$ and

\medskip

(c) $c^{{\bf C}^*}_{\rm top}(H^1({\bf P}^1,{\cal O}_{{\bf P}^1}(-d))
\otimes F^*\otimes E) =
\prod_{i=1}^e\prod_{j=1}^f\prod_{l=1}^{d-1}(x_i - y_j - l\hbar)$

\medskip

\nt Note that the {\it ratio} of $(b)$ and $(c)$ has many cancellations:

$$(b)/(c)
= (-1)^{ef(d-1)}\prod_{i,j}(y_j - x_i)(y_j - x_i + d\hbar)$$
Putting all of this together, we obtain:
$$e_{{\bf C}^*}(i_{\{d_i\}}Fl) = \frac{\prod_{i=1}^{k+1} \prod_{s=1}^{r_i}
\prod_{l=1}^{d_{m_{i}}} ( x_{m_{i-1}+s}+ l\hbar)^n} 
{\prod_{1 \le j < i \le k+1}
(-1)^{r_ir_j(d_{ij}-1)}\prod_{s,t}(x_{m_{i-1}+s} - 
x_{m_{j-1}+t} + d_{ij}\hbar)}$$
where $d_{ij} = d_{m_{i}} - d_{m_{j}}$ and $r_i = m_i - m_{i-1}$, 
and with Lemma 1.3:

$$J_d(\hbar) =
\sum_{\{d_i\}} \rho_*\left (\frac{\prod_{1 \le j < i \le k+1}
(-1)^{r_ir_j(d_{ij}-1)}
\prod_{s,t}(x_{m_{i-1}+s} - x_{m_{j-1}+t} + d_{ij}\hbar)}
{\prod_{i=1}^{k+1} \prod_{s=1}^{r_i}
\prod_{l=1}^{d_{m_{i}}} ( x_{m_{i-1}+s}+ l\hbar)^n}\right)$$
To complete the calculations, we need to 
compute the push-forwards.

\medskip

The following lemma is a special case of a general formula of Brion 
(\cite{Bri}):

\medskip

\nt {\bf Lemma 1.4:} Let $S$ be a rank $r$ vector bundle on a smooth variety
$X$, let $F:=Fl(m_1,...,m_k,S)$ be the associated flag bundle with projection
$\rho :F\rightarrow X$, and let $x_1,...,x_r$ denote 
the Chern roots of $S^*$ (and
of $\rho^*(S^*)$). Let $\rho_*$ be the push-forward map on Chow groups
(or rational (co)homology).
Then for any polynomial $P\in{\bf Q}[X_1,...,X_r]$,
$$\rho_*P(x_1,...,x_r)=\sum_{w}
w\left[ \frac {P(x_1,...,x_r)}{\prod_{1\leq j<i\leq k+1}\prod_{s,t}
(x_{m_{j-1}+t}-x_{m_{i-1}+s})}\right] ,$$
the sum over cosets $w\in S_r/(S_{r_1}\times ...
\times S_{r_{k+1}})$, where $S_r$ is the symmetric group, and 
the coset representatives are chosen to be the permutations 
$w\in S_r$ with descents in $\{ m_1,\dots ,m_k\}$.

\medskip

When we apply the Lemma to the formula for $J_d$, we obtain:
$$J_d=\sum_{\{d_i\}}
\sum_{w} w\left [{\displaystyle\frac{(-1)^{\sum_{j < i}r_ir_jd_{ij}}
\prod_{j < i}
\prod_{s,t}(x_{m_{i-1}+s} - x_{m_{j-1}+t} + d_{ij}\hbar)}
{\prod_{j<i}\prod_{s,t}
(x_{m_{i-1}+s}-x_{m_{j-1}+t})\prod_i \prod_s
\prod_{l=1}^{d_{m_{i}}} ( x_{m_{i-1}+s}+ l\hbar)^n}}\right ]$$

Notice first that the sign of each term in the
sum depends only on $r$ and $d$. Indeed,
$$\sum_{1\leq j<i\leq k+1}r_ir_j(d_{m_i}-d_{m_j}) 
\equiv (r-1)d\;\; ({\rm mod}\; 2)$$

\medskip

Second, we can simplify the formula for $J_d$ 
by replacing the double sum with a single sum over 
$r$-tuples $(d_1,\dots ,d_r)$ giving a partition
$d_1 + ... + d_r = d$. Namely, given such an $r$-tuple, let $k+1$
be the number of distinct $d_i$'s, with $r_1$ the multiplicity
of the smallest part, $r_2$ the multiplicity of the next smallest, etc.
Then there is
an unique $w\in S_r/(S_{r_1}\times ...\times S_{r_{k+1}})$ whose
inverse $w^{-1}$ arranges $(d_1,...,d_r)$ in nondecreasing order
$d_1\leq d_2\leq ...\leq d_r$ and we have:
$$\frac{\prod_{1\leq j<i\leq r}(x_i-x_j+(d_i-d_j)\hbar)}
{\prod_{1\leq j<i\leq r}(x_i-x_j)\prod_{i=1}^r\prod_{l=1}^{d_i}
(x_i+l\hbar)^n}=$$
$$w\left [ \frac{\prod_{1\leq j < i\leq k+1}
\prod_{s,t}(x_{m_{i-1}+s} - x_{m_{j-1}+t} + d_{ij}\hbar)}
{\prod_{1\leq j<i\leq k+1}\prod_{s,t}
(x_{m_{i-1}+s}-x_{m_{j-1}+t})\prod_{i=1}^{k+1} \prod_{s=1}^{r_i}
\prod_{l=1}^{d_{m_i}} ( x_{m_{i-1}+s}+ l\hbar)^n} \right ] ,$$
so that, putting all this together, we arrive at:

\medskip

\nt {\bf Theorem 1.5:} The $J$-function of the Grassmannian 
${\bf G} = G(r,n)$ is
$$J^{\bf G} = e^{\frac {t\sigma_1}{\hbar}}
\sum_{d\geq 0}e^{dt}J_d(\hbar), \ \mbox{where}$$
$$J_d(\hbar)=(-1)^{(r-1)d}\sum_{\stackrel{(d_1,\dots ,d_r)}{d_1+...+d_r=d}}
\frac{\prod_{1\leq i<j\leq r}(x_i-x_j+(d_i-d_j)\hbar)}
{\prod_{1\leq i<j\leq r}(x_i-x_j)\prod_{i=1}^r\prod_{l=1}^{d_i}
(x_i+l\hbar)^n} ,$$
and $x_1,...x_r$ are the Chern roots of $S^*$, the dual of the
tautological subbundle.

\medskip

\nt {\bf Remark:} Theorem 1.5 specializes to Givental's formula for
${\bf P}^{n-1}\cong G(1,n)$:
$$J^{{\bf P}^{n-1}} = e^{\frac {tx}{\hbar}}\sum_{d \ge 0} \frac {e^{dt}}
{\prod_{l = 1}^d(x + l\hbar)^n}$$ 

\medskip

{\bf First Proof of the Hori-Vafa Conjecture:} Simply compute:
$$J^{\bf P} = \prod_{i=1}^r J^{{\bf P}^{n-1}} = 
e^{\frac{t_1x_1 + ... + t_rx_r}{\hbar}} \sum_{(d_1,...,d_r)} 
\frac{e^{d_1t_1 + ... + d_rt_r}}{\prod_{i=1}^r
\prod_{l=1}^{d_i}(x_i + l\hbar)^n}$$
so
$$\frac{{\cal D}_{\Delta}(J^{\bf P})}{\Delta} = 
e^{\frac {t_1x_1 + ... + t_rx_r}{\hbar}}\sum_{(d_1,...,d_r)}
\frac{e^{d_1t_1 + ... + d_rt_r}\prod_{1\leq i<j\leq r}(x_i-x_j+(d_i-d_j)\hbar)}
{\prod_{1\leq i<j\leq r}(x_i-x_j)\prod_{i=1}^r\prod_{l=1}^{d_i}
(x_i+l\hbar)^n} ,$$
and then the conjecture immediately follows. 

\medskip

Finally, there is an analogous formula for the 
$J$-function in $T$-equivariant cohomology (the action of 
$T = ({\bf C}^*)^n$ was described in 
the proof Lemma 1.3):  
$$J_d^T(\hbar,\lambda_1,\dots ,\lambda_n)=
{ev_T}_*\frac 1{e_{{\bf C}^*\times T}(F_0)}
\in H_T^*(G(r,n),{\bf Q}).$$ 

Indeed, a computation analogous to the proof of Theorem 1.5 gives:

\medskip

\nt {\bf Theorem 1.5':} The equivariant $J$-function of the Grassmannian 
$G(r,n)$ is
$${J^{\bf G}}^T = e^{\frac {t\sigma_1}\hbar}\sum_{d\geq 0}e^{dt}J_d^T ,$$
where
$$J_d^T=(-1)^{(r-1)d}\sum_{\stackrel{(d_1,\dots ,d_r)}{d_1+...+d_r=d}}
\frac{\prod_{1\leq i<j\leq r}(x_i-x_j+(d_i-d_j)\hbar)}
{\prod_{1\leq i<j\leq r}(x_i-x_j)\prod_{i=1}^r\prod_{l=1}^{d_i}
\prod _{j=1}^n (x_i-\lambda _j+ l\hbar )} ,$$

The substitution of $\prod_{l=1}^{d_i}
\prod _{j=1}^n (x_i-\lambda _j+ l\hbar )$ in place of $\prod_{l=1}^{d_i}
(x_i+l\hbar )^n$ corresponds to taking into
account the $T$-action on the trivial bundles ${\bf C}^n$ in the
Euler sequences. Under this action, the trivial bundle decomposes into
a direct sum of representations ${\bf C}^n=\bigoplus_{i=1}^n{\bf C}_i $,
with $c_1^T({\bf C}_i)=-\lambda_i$. Consequently, 
the $T$-equivariant version of formula $(**)$ is
$$ c_{\rm top}(T{\rm Quot}/
TFl) = \frac{(\prod _{i=1}^n c_{\rm top}(\pi_*{(\cal K}^*/\pi^*S^* ) 
\otimes {\bf C}_i )) \cdot
c_{\rm top}(R^1\pi_*({\cal K}^*\otimes {\cal K}))}
{c_{\rm top}(\pi_*({\cal K}^*\otimes
{\cal K})/K)},$$
where the Chern classes are now ${\bf C}^*\times T$-equivariant.

\newpage

\nt {\bf \S 2. A Quantum Cohomology Proof of the Hori-Vafa Conjecture.} 
This proof follows from
``quantum''  versions of some results of \cite{ES} 
and \cite{Mar}, on 
the classical cohomology.
Consider the rational map:
$$\Phi: {\bf P} = \prod_{i=1}^r {\bf P}^{n-1} --\!\!\!> {\bf G}; 
(p_1,...,p_r) \mapsto
\mbox{span}\{p_1,...,p_r\}
\subset {\bf C}^n$$ 

Recall the the presentation of the cohomology ring of ${\bf G}$: 
$$H^*({\bf G})\cong {\bf C}[x_1,\dots ,x_r]^{S_r}/\langle h_{n-r+1},
\dots h_n\rangle,$$
and the presentation of the cohomology ring of ${\bf P}$:
$$H^*({\bf P})\cong {\bf C}[x_1,\dots ,x_r]/\langle x_1^n,\dots ,x_r^n
\rangle.$$
Thus each class $\gamma\in
H^*({\bf G})$ together with a symmetric polynomial $P(x_1,...,x_r)$ 
representing $\gamma$ lifts to
${\tilde{\gamma}}\in H^*({\bf P})$ by evaluating $P$ in $H^*({\bf P})$. 

\medskip

For a partition 
$\mu=(n-r\geq \mu_1\geq\dots\geq\mu_r\geq 0)$ one defines
$\sigma_{\mu}$ to be the corresponding Schur
polynomial:
$$\sigma_\mu = \frac{\mbox{det}\left(
\begin{array}{cccc}x_1^{\mu_1 + r-1} & x_1^{\mu_2 + r - 2} & 
\cdots & x_1^{\mu_r} \\
x_2^{\mu_1 + r-1} & x_2^{\mu_2 + r - 2} & \cdots & x_2^{\mu_r} \\
& & \vdots \\ 
x_r^{\mu_1 + r-1} & x_r^{\mu_2 + r - 2} & \cdots & x_r^{\mu_r}
\end{array}\right)}
{\mbox{det}\left(
\begin{array}{cccc} x_1^{r-1} & x_1^{r - 2} & \cdots & 1 \\
x_2^{r-1} & x_2^{r - 2} & \cdots & 1 \\
& & \vdots \\ 
x_r^{r-1} & x_r^{r - 2} & \cdots & 1
\end{array}\right)}$$ 

In the 
cohomology of ${\bf G}$ this represents a ``Schubert'' class, and such
classes form an additive 
basis of $H^*({\bf G})$ as $\mu$ runs over the set of partitions.
Since each $x_i$ appears in the Schur polynomial $\sigma_{\mu}$
with exponent at most $n-r$,
the lift of a Schubert class to $H^*({\bf P})$
determines the Schur polynomial. We will
often identify a partition with its Young diagram. The partitions
corresponding to Schubert classes on ${\bf G}$ are those whose Young diagrams
fit in an $r\times (n-r)$ rectangle.
Finally, recall that the denominator of $\sigma_\mu$ is the Vandermonde 
determinant:
$$\Delta = \prod_{1\leq i< j\leq r}(x_i-x_j)$$

\medskip

The relation between the cohomology of ${\bf G}$ and that of ${\bf P}$
is encoded in the following ``integration formula''

\medskip

\nt {\bf Theorem 2.1 (Martin):} For any cohomology
class $\gamma\in H^*({\bf G})$,
$$\int_{\bf G} \gamma=\frac{(-1)^{r\choose 2}}{r!}\int_{{\bf P}}
{\tilde{\gamma}}\cup_{{\bf P}}\Delta^2.$$

\medskip

\nt {\bf Corollary 2.2 (Ellingsrud-Str\o mme):} The linear map
$$\theta :H^*({\bf G})\longrightarrow
H^*({\bf P}); \ \theta(\gamma)={\tilde{\gamma}}\cup_{{\bf P}}\Delta$$ 
is injective, and its image $V \subset H^*({\bf P})$ is the subspace of 
anti-symmetric classes.

\medskip

{\bf Proof:} If $\theta(\gamma)=0$, then $$\int_{{\bf P}}
({\tilde{\gamma}}\cup_{{\bf P}}\Delta)\cup_{{\bf P}} b=0,\; 
\forall\; b\in H^*({\bf P}).$$
In particular $$\int_{\bf G} \gamma\cup_{\bf G}\gamma '=\int_{{\bf P}}
({\tilde{\gamma}}\cup_{{\bf P}}\Delta)\cup_{{\bf P}}
({\tilde{\gamma '}}\cup_{{\bf P}}\Delta)=0,\;
\forall\; \gamma '\in H^*({\bf G}),$$
by the integration formula. Hence $\gamma =0$. Finally, it is evident that
each class $\tilde \gamma \cup_{\bf P} \Delta$ in the image of $\theta$ is 
anti-symmetric, but conversely, an anti-symmetric class is always of 
the form $\tilde \gamma \cup_{\bf P} \Delta$, with $\tilde \gamma \in {\bf
C}[x_1,...,x_r]^{S_r}$.

\medskip

Thus $V={\rm span}\left( \{\sigma_{\mu}\cup_{\bf P}\Delta\mid \mu=(n-r\geq\mu_1\geq\dots
\geq\mu_r\geq 0)\}\right )$. In addition, let $W = V^\perp$ be its
orthogonal complement with respect to the intersection form.
Note that the polynomial $\sigma_{\mu}\Delta$ represents the class
$\sigma_{\mu}\cup_{{\bf P}}\Delta$.

\medskip

If $\check{\mu}$ is the dual partition
(the complement of $\mu$ in the $r\times(n-r)$ rectangle),
then it follows
from the integration formula that
$$\int_{{\bf P}}(\sigma_{\mu}\Delta)\cup(\sigma_{\nu}\Delta)=
(-1)^{r\choose 2}r!\delta_{\nu\check{\mu}}$$
Choose a basis $\{ b_i\}$ of $W$ that extends $\{\sigma_\mu\Delta\}$ to a basis ${\cal B}$ of 
$H^*({\bf P})$ with intersection form $((-1)^{r\choose 2}r!){\rm Id}$ and let
$$p_V:H^*({\bf P})\longrightarrow V$$ be the orthogonal projection onto $V$.

\medskip

\nt {\bf Corollary 2.3.} For any two partitions $\mu, \nu$,
$$\theta(\sigma_{\mu}\cup_{\bf G}\sigma_{\nu})=\theta(\sigma_{\mu})
\cup_{{\bf P}}\sigma_{\nu}.$$

{\bf Proof:} Since the right hand side is antisymmetric, it lies
in $V$, hence
$$\theta(\sigma_{\mu})\cup_{{\bf P}}\sigma_{\nu}=\sum_{\rho}
\left(\int_{{\bf P}}(\sigma_{\mu}\Delta)\cup\sigma_{\nu}
\cup(\sigma_{\rho}\Delta)
\right)\frac{(-1)^{r\choose 2}}{r!} \sigma_{\check{\rho}}\Delta $$
so continuing with Martin's formula gives:
$$\theta(\sigma_{\mu})\cup_{{\bf P}}\sigma_{\nu}=
\sum_{\rho}\left(\int_{\bf G}\sigma_{\mu}
\cup\sigma_{\nu}\cup\sigma_{\rho}\right)\theta( \sigma_{\check{\rho}})=
\theta(\sigma_{\mu}\cup_{\bf G}\sigma_{\nu}).$$

This corollary has a natural generalization to quantum cohomology.

\medskip 

Recall that the 
3-point genus zero Gromov-Witten invariants of  a smooth projective 
variety $X$ 
define the small quantum cohomology ring as a deformation
of the usual cup-product on $X$. Specifically,
if $\{\gamma_i\}_{i=1}^s$ is a ${\bf C}$-basis of $H^{2*}(X)$ and 
$\{\check \gamma_i\}_{i=1}^s$ is the dual basis (with respect to the 
intersection form) then the quantum product of two cohomology classes is:
$$\gamma * \delta =\sum_{\beta\in H_2(X,{\bf Z})} \sum_k e^{\beta}
\langle\gamma,\delta,\gamma_k\rangle^X_{\beta} \ {\check \gamma_k},$$

We think of $\beta$ either as the class of an algebraic curve, and compute:
$$\langle\gamma,\gamma',\gamma ''\rangle^X_{\beta} = 
\int_{[\overline M_{0,3}(X,\beta)]} 
ev_1^*(\gamma) \cup ev_2^*(\gamma')\cup ev_3^*(\gamma'')$$
where 
$$ev_i: \overline M_{0,n}(X,\beta) \rightarrow X;\ ev([f]) = f(p_i)$$
and $[\overline M_{0,3}(X,\beta)] \in A_*(\overline M_{0,3}(X,\beta))_{\bf Q}$
is the ``virtual fundamental class'' in the Chow group of the 
Kontsevich-Manin stack of stable maps (see \cite{BF}). Otherwise, we
think of
$\beta =\sum_{i=1}^m d_it_i$ as an  element of the dual space to 
$H^2(X,{\bf C})$ so that:
$$e^\beta = e^{d_1t_1 + ... + d_mt_m}$$

If we assume 
that $H^2(X)$ has a basis 
consisting of nef divisors, this 
deformation defines small quantum cohomology as an algebra 
over ${\bf C}[[e^{t_1},...,e^{t_m}]]$ ($m = \mbox{dim}(H^2(X))$). In the cases
$X={\bf G}$ and
$X={\bf P}$, set:
$$QH^*({\bf G}) := H^*({\bf G})\otimes_{{\bf C}} {\bf C}[e^t]$$
with multiplication defined on the Schubert basis by
$$\sigma_{\mu}*_{\bf G}\sigma_{\nu}=
\sum_{d\geq 0}\sum_{\rho}
e^{dt}\langle \sigma_{\mu},\sigma_{\nu},\sigma_{\rho}\rangle^{\bf G}_d
\ \sigma_{{\check\rho}}$$ 
and similarly:
$$QH^*({\bf P})= H^*({\bf P})\otimes_{{\bf C}}{\bf C}[e^{t_1},...,e^{t_r}]$$
with multiplication given on a basis $\{\gamma_K\}$ of $H^*({\bf P})$ by
$$\gamma_I*_{\bf P}\gamma_J=\sum_{d_1,\dots ,d_r\geq 0}\sum_K
e^{d_1t_1 + ... + d_rt_r} \langle\gamma_I,\gamma_J,\gamma_K
\rangle^{\bf P}_{(d_1,...,d_r)}{\check
\gamma_K}.$$  where $(d_1,...,d_r)$ is the multi-degree of a curve in 
${\bf P} = ({\bf P}^{n-1})^r$.


We've already mentioned the well-known presentations:
$$QH^*({\bf G})\cong{\bf C}[x_1,\dots ,x_r]^{S_r}[e^t]/
\langle h_{n-r+1},\dots h_n - (-1)^{r-1}e^t\rangle,$$
and:
$$QH^*({\bf P}) \cong {\bf C}[x_1,\dots ,x_r] [e^{t_1},\dots
,e^{t_r}]/\langle x_1^n-e^{t_1},\dots , x_r^n-e^{t_r}\rangle.$$

We will use two bases 
for $QH^*({\bf P})$ (over ${\bf C}[e^{t_1},...,e^{t_r}]$). One is the 
basis ${\cal B}$ above. The other is the natural
one consisting of monomials
$X^I=x_1^{i_1}\dots x_r^{i_r}$, $0\leq r_j\leq n-1$. In the latter
basis, the quantum product $X^I*X^J$
is obtained simply by multiplying $X^I\cdot X^J
= X^{I+J}$
and then substituting: 
$$x_i^n = e^{t_i}$$
As a consequence of this we have

\medskip

\nt {\bf Lemma 2.4:} For any partition 
$\mu=(n-r \ge \mu_1 \ge \mu_2 \ge ... \ge \mu_r \ge 0)$
$$\sigma_{\mu}*_{{\bf P}}\Delta=\sigma_{\mu}\cup_{{\bf P}}\Delta=
\sigma_{\mu}\Delta.$$

{\bf Proof:} The highest power of each $x_i$ that appears
in $\sigma_{\mu}\Delta$ is at most $n-1$, as is evident 
from the Schur polynomial.

\medskip

The action of the symmetric group $S_r$ on polynomials
in the $x$'s and $e^t$'s gives an action on $QH^*({\bf P})$.
Let $\overline{QH}^*({\bf P})$ denote the quotient obtained by substituting
$$e^{t_i}=(-1)^{r-1}e^t;\ \ i=1,\dots ,r$$

The $S_r$-action descends to an action on
$\overline{QH}^*({\bf P})$ (permuting only the $x_i$'s) and the anti-invariants
are $V\otimes {\bf C}[e^t]$. Now extend $\theta$ by linearity
over ${\bf C}[e^t]$ to:
$$\overline\theta: QH^*(G)\longrightarrow \overline{QH}^*({\bf P}),$$
which obviously continues to be injective with image equal to 
$V\otimes {\bf C}[e^t]$. 

\medskip
The quantum analogue of Corollary 2.3 is the main result of this section:

\medskip

\nt {\bf Theorem 2.5:} For any two partitions $\mu, \nu$:
$$
\overline\theta(\sigma_{\mu}*_{\bf G}\sigma_{\nu})=(\theta(\sigma_{\mu})
*_{{\bf P}}\sigma_{\nu})|_{e^{t_i}=(-1)^{r-1}e^t}.$$
In the spirit of this paper we will give two proofs. 

\medskip

{\bf First Proof:} For a sequence of non-negative integers
$\alpha=(\alpha_1,\alpha_2,\dots ,\alpha_r)$ define
$D_{\alpha}={\rm det}(x_{i}^{\alpha_{j}})_{1\leq i,j\leq r}$ so 
if we let $\delta = (r-1,r-2,....,0)$, then:
$$D_{\delta} = \Delta, \ \sigma_\mu(x_1,...,x_r) = 
\frac{D_{\mu + \delta}}{D_\delta} \
\mbox{for all partitions $\mu$,}$$
and in general $D_\alpha \ne 0$ if and only if the 
$\alpha_j$'s are all distinct.

\medskip

The product of two Schur polynomials in $r$ variables is given by
\begin{equation}\sigma_{\mu}\sigma_{\nu}=
\sum_{\rho}N_{\mu,\nu}^{\rho}\sigma_{\rho},\label{LR}\end{equation}
where $N_{\mu,\nu}^{\rho}$ are the Littlewood-Richardson coefficients and 
$\rho$ is constrained to be of the form 
$\rho_1 \ge \rho_2 \ge ... \ge \rho_r \ge 0$,
but we do not require $\rho_1 \leq n-r$.  Thus:
\begin{equation}\sigma_{\mu}\Delta\sigma_{\nu}=
\sum_{\rho}N_{\mu,\nu}^{\rho}D_{\rho+\delta}.\label{projLR}\end{equation} 

Let $$I_{\bf G} =\langle h_{r+1},\dots ,h_n-(-1)^{r-1}e^t\rangle$$ and 
$$I_{{\bf P}}=\langle x_1^n-(-1)^{r-1}e^t,\dots ,x_r^n-(-1)^{r-1}e^t\rangle $$ 
be the ideals of relations from
the presentations of $QH^*(G)$ and of $\overline{QH}^*({\bf P})$. 
Then $\overline{\theta}(\sigma_{\mu} *_{\bf G}\sigma_{\nu})$ is obtained by
reducing the right-hand side of (\ref{LR}) modulo $I_{\bf G}$ and multiplying
the result by $\Delta$, while $\theta(\sigma_{\mu})*_{{\bf P}}\sigma_{\nu}$
is obtained by reducing the right-hand side of (\ref{projLR}) modulo
$I_{{\bf P}}$.

\medskip

The reduction of a Schur polynomial modulo $I_G$ can be done
using the rim-hook algorithm of \cite{BCF}, which goes as follows.
The rim of a partition $\rho=(\rho_1,\dots ,\rho_r)$ consists of the boxes
on the border of its Young diagram from 
the northeast (upper-right) to the southwest (lower-left) corners. 
An $n$-rim hook
is a contiguous collection of $n$-boxes in the rim 
with the property that when removed from the Young
diagram of $\rho$ the remaining shape
is again the diagram of a partition. 
If we start at the end of any row and count $n$
boxes along the rim, moving down and to the
left, then we either obtain an $n$-rim hook, or
the process ends in a box directly above the last box
in some row. In the latter case, the set of $n$ boxes
is called an illegal $n$-rim. With these definitions,
the Main Lemma from \cite{BCF} gives:

\medskip

$(i)$ If $\rho_1>n-r$ and either $\rho$ contains no $n$-rim hook,
or $\rho$ contains an illegal $n$-rim, then: 
$$\sigma_{\rho} \equiv 0 \ \ ({\rm mod}\ \  I_G)$$

\medskip

$(ii)$ If $\pi=(\pi_1\geq\pi_2\geq\dots\geq\pi_r\geq 0)$
is obtained from $\rho$ by removing an $n$-rim hook,
then 
$$\sigma_{\rho} \equiv (-1)^{r-s}e^t\sigma_{\pi}\ \ ({\rm mod}\ \ I_G)$$
where $s$ is the height of the rim hook, i.e., the number
of rows it occupies.

\medskip

On the other hand, it is obvious that if $\alpha_j\geq n$,
then
$$D_{\alpha} \equiv (-1)^{r-1}e^t D_{\alpha '}\ \ ({\rm mod}\ \ 
I_{{\bf P}}),$$
where $\alpha '=(\alpha_1,\dots ,\alpha_j-n,\dots ,\alpha_r)$.
It is easy to see that if $\rho_1>n-r$ and $\rho$
contains no $n$-rim, then $\rho +\delta -(n,0,\dots ,0)$ has
two equal entries, hence $D_{\rho+\delta}\equiv 0 \ ({\rm mod}\ I_{{\bf P}})$.
Assume now that $\rho$ has an $n$-rim, starting in row $j$ and
ending in row $k$ for some $1\leq j<k\leq r$ (so that its height
is $s:=k-j+1$).
Then we must have $\rho_j+\delta_j>n$. Moreover, if the rim is
illegal, then $k<r$ and
$$\rho_j=\rho_{k+1}+n-(k+1)+j.$$
This implies that 
$$\rho_j+\delta_j-n=\rho_{k+1}+\delta_{k+1},$$
so again $D_{\rho+\delta} \equiv 0 \ ({\rm mod}\ I_{{\bf P}})$.

\medskip

Finally, if the $n$-rim is a rim hook, let $\pi$ be the partition
obtained by removing it from $\rho$.
Then 
one gets $\pi+\delta$ from 
$\rho+\delta-(0,\dots ,0,n,0,\dots ,0)$ by taking its $j$th entry
$\rho_j+\delta_j-n$ and moving it past the next $s-1$ entries. This gives
$$D_{\rho+\delta}=(-1)^{s-1}(-1)^{r-1}e^tD_{\pi+\delta} \equiv
(-1)^{r-s}e^tD_{\pi + \delta}  \ ({\rm mod} \ I_{{\bf P}}),$$
and finishes the first proof of the Theorem. 

\medskip

{\bf Second Proof of Theorem 2.5:} In the basis 
${\cal B} = \{\sigma_\rho\Delta\} \cup\{b_i\}:$

$$\theta(\sigma_{\mu})
*_{{\bf P}}\sigma_{\nu} = \frac{(-1)^{r \choose 2}}{r!}
\left(\sum_{\rho,\{ d_i\}}e^{d_1t_1 + ... + d_rt_r}
\langle\sigma_{\mu}\Delta,
\sigma_{\nu},\sigma_{\rho}\Delta\rangle^{\bf P}_{\{d_i\}}
\sigma_{\check{\rho}}\Delta\right.$$ 
$$+ \left.\sum_{j,\{ d_i\}}e^{d_1t_1 + ... + d_rt_r}
\langle\sigma_{\mu}\Delta,
\sigma_{\nu},b_j\rangle^{\bf P}_{\{d_i\} }
{\check{b_j}}\right)$$

After substituting for the $e^{t_i}$'s this becomes

$$\frac{(-1)^{{r \choose 2}}}{r!}\left(
\sum_{\rho ,d}\sum_{d_1+\dots +d_r=d}
(-1)^{d(r-1)}e^{dt}\langle \sigma_{\mu}\Delta,
\sigma_{\nu},\sigma_{\rho}\Delta\rangle^{\bf P}_{\{d_i\} }
\sigma_{\check{\rho}}\Delta
\right.$$

$$\left. +\sum_{j,d}\sum_{d_1+\dots +d_r=d}
(-1)^{d(r-1)}e^{dt}\langle \sigma_{\mu}\Delta,\sigma_{\nu},b_j
\rangle^{\bf P}_{\{d_i\} }
{\check{b_j}}\right ).$$

This is antisymmetric in the $x_i$'s,
so each:
$$\sum_{d_1+\dots +d_r=d}
\langle \sigma_{\mu}\Delta,\sigma_{\nu},b_j
\rangle^{\bf P}_{\{d_i\} }=0$$

Moreover, equating the coefficients for
the basis elements
$\sigma_{\check{\rho}}\Delta$ gives us the following equivalent formulation of 
Theorem 2.5:
\begin{equation}
\langle \sigma_{\mu},\sigma_{\nu},\sigma_{\rho}\rangle^{\bf G}_d =
\frac{(-1)^{r \choose 2}}{r!}\sum_{d_1+\dots +d_r=d}
(-1)^{d(r-1)}\langle \sigma_{\mu}\Delta,
\sigma_{\nu},\sigma_{\rho}\Delta\rangle^{\bf P}_{\{d_i\} }
\label{qintegration}\end{equation}
which is a ``quantum'' analogue of the integration formula (Theorem 2.1).

\medskip

This can be proved directly with the 
Vafa - Intriligator residue formula. Given classes 
$\sigma_\mu,\sigma_\nu,\sigma_\rho$ with $|\mu| + |\nu| + |\rho| = 
\mbox{dim}({\bf G})$, 
then the ``classical'' version of the formula gives:
$$\int_{\bf G} \sigma_\mu \sigma_\nu \sigma_\rho = 
\frac{(-1)^{\left( r \atop 2 \right)}}{r!n^r} 
\sum_{\{(\epsilon_1,...,\epsilon_r) | \epsilon_i^n = 1\}}
\sigma_\mu\sigma_\nu \sigma_\rho \prod_i \epsilon_i \prod_{i \ne j} 
(\epsilon_i - \epsilon_j)$$
where $\sigma_\mu = \sigma_\mu(\epsilon_1,...,\epsilon_r)$ on the right side,
and likewise for $\sigma_\nu, \sigma_\rho$. 
The quantum formula is essentially
identical. If $|\mu| + |\nu| + |\rho| = nd + \mbox{dim}({\bf G})$ then:
$$\langle \sigma_\mu, \sigma_\nu,\sigma_\rho\rangle ^{\bf G}_d = 
\frac{(-1)^{\left( r \atop 2 \right) + d(r-1)}}{r!n^r} 
\sum_{\{(\epsilon_1,...,\epsilon_r) | \epsilon_i^n =
1\}}\sigma_\mu\sigma_\nu \sigma_\rho \prod_i \epsilon_i \prod_{i \ne j} 
(\epsilon_i - \epsilon_j)$$

This formula is proved (see \cite{Ber1}) by taking the potential function:
$$P = \sum_{i=1}^r \left( \frac{{x_i}^{n+1}}{n+1} - (-1)^{r-1}e^tx_i\right)$$
satisfying $\langle \frac{\partial P}{\partial \sigma_1},...,
\frac{\partial P}{\partial
\sigma_n}\rangle = I_{\bf G}$ and summing over the residues of the $r$-form:
$$\sigma_\mu \sigma_\nu \sigma_\rho \frac{d\sigma_1 
\wedge ... \wedge d\sigma_r}
{\prod_i \frac {\partial P}{\partial \sigma_i}}$$ 

On the other hand, this same $P$ satisfies 
$\langle \frac{\partial P}{\partial x_1},...,\frac{\partial P}
{\partial x_n}\rangle =
I_{\bf P}$, and it follows from Lemma 2.4 and summing over 
residues  of the $r$-form:
$$\sigma_\mu \Delta \sigma_\nu \sigma_\rho\Delta 
\frac{dx_1 \wedge ... \wedge dx_r}
{\prod_i \frac {\partial P}{\partial x_i}}$$
(just as in \cite{Ber1}) that:
$$\sum_{d_1 + ... + d_r = d} 
\langle \sigma_\mu\Delta, \sigma_\nu, \sigma_\rho \Delta
\rangle_{\{d_i\}}^{\bf P} = 
\frac{(-1)^{\left( r \atop 2 \right)}}{n^r} 
\sum_{\{(\epsilon_1,...,\epsilon_r) | \epsilon_i^n = 1\}}
\sigma_\mu \Delta \sigma_\nu \sigma_\rho \Delta \prod_{i} \epsilon_i$$
and this, plus $\Delta^2 = \prod_{i \ne j} (x_i - x_j)$, 
proves the Theorem.

\medskip

We next want to consider the flat sections for ${\bf P}$, that is, 
solutions to:
$$\hbar\frac{\partial}{\partial t_i}f=x_i*_{{\bf P}}f,\;\;\; i=1,\dots ,r.$$
with
$$f(t_1,\dots t_r, \hbar) = 
\sum_{\mu} f_\mu (\sigma_{\mu}\Delta) + \sum_j f_jb_j$$ 
a vector field (over $H^2({\bf P})$) written in terms of 
the basis ${\cal B}$. (Here, and for the
rest of the paper, we abuse the notation slightly by writing vector fields
as cohomology-valued, via the canonical identification of the tangent space
at a point in $H^*(X)$ with $H^*(X)$ itself.) Let:
$$p_V(f) = \sum_{\mu} f_\mu (\sigma_{\mu}\Delta)\ \mbox{with}\ 
p_V(f)^\perp = \sum_j f_jb_j$$
for arbitrary vector fields. Then:

\medskip

\nt {\bf Lemma 2.6:} For all 
$\nu = (n-r \ge \nu_1 \ge ... \ge \nu_r \ge 0)$ 
and all  $f = f(t_1,...,t_r,\hbar)$:
$$p_V\left( \left (\sigma_{\nu}*_{{\bf P}}f\right )
|_{t_i=t+\pi\sqrt{-1}(r-1)}\right)
=\left( \sigma_{\nu}*_{{\bf P}}p_V(f)\right)|_{t_i=t+\pi\sqrt{-1}(r-1)}$$

{\bf Proof:} We need to show that:
$$p_V\left( \left (\sigma_{\nu}*_{{\bf P}}p_V(f)^{\perp}\right )
|_{t_i=t+\pi\sqrt{-1}(r-1)}\right)=0.$$

But
$\sigma_{\nu}*_{{\bf P}}p_V(f)^{\perp}=\sum_jf_j(\sigma_{\nu}
*_{{\bf P}}b_j)$  and after specializing the $e^{t_i}$'s,
the coefficient of $\sigma_{\rho}\Delta$ in $\sigma_{\nu}
*_{{\bf P}}b_j$ is

$$\sum_{d_1+\dots +d_r=d}
\langle \sigma_{\nu},b_j,\sigma_{\check{\rho}}\Delta
\rangle_{\{d_i\} }^{\bf P}$$
which vanishes, as noted earlier.

\medskip

Now let $f$ be a flat section for ${\bf P}$. Then:
$$\delta f=\sigma_1*_{{\bf P}}f \ \mbox{where} \ \delta =
\hbar\sum_{i=1}^r \frac{\partial}{\partial t_i} \ \mbox{and}\ 
\sigma_1= x_1+\dots +x_r$$
Since $p_V(\delta f)=\delta p_V(f)$, we get 
$\delta p_V(f)=p_V(\sigma_1*_{{\bf P}}f)$ and:
$$(\delta p_V(f))|_{t_i=t+\pi\sqrt{-1}(r-1)}=
(\sigma_1*_{{\bf P}}p_V(f))|_{t_i=t+\pi\sqrt{-1}(r-1)}$$
from Lemma 2.6 applied to $\nu = (1 \ge 0 \ge 0 ... \ge 0)$.

\newpage

\nt {\bf Theorem 2.7:} If $f=\sum_{\mu}f_{\mu}(\sigma_{\mu}\Delta)+
\sum_jf_jb_j$ is a flat section for ${\bf P}$, then

$$\left. \frac{p_V(f)}{\Delta}\right |_{t_i=t+\pi\sqrt{-1}(r-1)}
=\sum_{\mu}\left( f_{\mu}|_{t_i=t+\pi\sqrt{-1}(r-1)}\right)\sigma_{\mu}$$
is a flat section for ${\bf G}$.

\medskip

{\bf Proof:} We need to show that
$$\hbar\frac{d}{dt}\left( \left. \frac{p_V(f)}{\Delta}
\right |_{t_i=t+\pi\sqrt{-1}(r-1)}\right)-\sigma_1*_{\bf G}
\left( \left. \frac{p_V(f)}{\Delta}
\right |_{t_i=t+\pi\sqrt{-1}(r-1)}\right)=0.$$
Think of this as an $H^*({\bf G})$-valued function. By Corollary 2.2, 
it suffices to show

$$\overline{\theta}\left(\hbar\frac{d}{dt}\left( \left. \frac{p_V(f)}{\Delta}
\right |_{t_i=t+\pi\sqrt{-1}(r-1)}\right)-\sigma_1*_{\bf G}
\left( \left. \frac{p_V(f)}{\Delta}
\right |_{t_i=t+\pi\sqrt{-1}(r-1)}\right)\right)=0.$$

Since $\overline{\theta}$ commutes with $\hbar(d/dt)$, 
Theorem 2.5 allows us to rewrite this as:

$$\hbar\frac{d}{dt}\left(p_V(f)|_{t_i=t+\pi\sqrt{-1}(r-1)}\right)-
\left(\left. \sigma_1\Delta *_{{\bf P}}\frac{p_V(f)}{\Delta}\right)
\right |_{t_i=t+\pi\sqrt{-1}(r-1)}=0$$

But now by Lemma 2.4:
$$\sigma_1\Delta=\sigma_1*_{{\bf P}}\Delta\ \mbox{and}$$
$$\left(\Delta *_{{\bf P}}
\left. \frac{p_V(f)}{\Delta}
\right)\right |_{t_i=t+\pi\sqrt{-1}(r-1)}=
p_V(f)|_{t_i=t+\pi\sqrt{-1}(r-1)}.$$
so that we may rewrite the desired relation once again 
(using the chain rule) as:
$$(\delta p_V(f))
|_{t_i=t+\pi\sqrt{-1}(r-1)} - 
(\sigma_1*_{\bf P}p_V(f))|_{t_i=t+\pi\sqrt{-1}(r-1)} = 0$$
and the validity of this relation, as noted above, follows from Lemma 2.6.

\medskip

Finally, we are ready for the:

\medskip

\nt {\bf Quantum Cohomology Proof of the Hori-Vafa Conjecture:}

\medskip

The conjecture will follow from Theorem 2.7 and from basic properties
of the explicit fundamental matrices of flat sections for ${\bf P}$ and 
${\bf G}$
given by two-point gravitational descendants, which we now recall.

For a smooth projective $X$, 
the two-point invariants involved will be of the form

$$\langle\psi^m\gamma,\gamma'\rangle^X_{\beta} = 
\int_{[\overline M_{0,2}(X,\beta)]} 
\psi^m ev_1^*(\gamma) \cup ev_2^*(\gamma'),$$
where $\gamma ,\gamma'$ are (even) cohomology classes on $X$
and $\psi$ is the first Chern class of the line bundle $L_1$ on
$\overline M_{0,2}(X,\beta)$ with fiber over the point
$[C,p_1,p_2,f]$ the cotangent line to $C$ at $p_1$.
More formally, $L_1=s_1^*(\omega_{\pi})$, with $s_1$ the natural section --
corresponding to the first marked point -- of the projection
$\pi$ from the universal curve to $\overline M_{0,2}(X,\beta)$, and 
$\omega_\pi$ the relative dualizing sheaf.
If, as before, $\{T_i\}_{i=0}^s$ and 
$\{{\check T_i}\}_{i=0}^s$ are dual bases of $H^{2*}(X)$,
with $\{T_1,\dots ,T_m\}$ a basis of $H^2(X)$, then for
each $\gamma\in H^{2*}(X)$, the cohomology-valued
function 
$$F_{\gamma}^X(t,e^t,\hbar):=\sum_{k=0}^s\langle\langle
\frac{\gamma e^{(t_1T_1+\dots t_mT_m)/\hbar}}{\hbar-\psi},T_k
\rangle\rangle^X{\check T_k}$$
is a flat section, where we used the customary notation 
$$\langle\langle\frac{\gamma e^{(t_1T_1+\dots t_mT_m)/\hbar}}{\hbar-\psi},T_k
\rangle\rangle^X=\sum_{\beta}e^{\beta}\langle
\frac{\gamma e^{(t_1T_1+\dots t_mT_m)/\hbar}}{\hbar-\psi},T_k
\rangle_{\beta}^X.$$
When $\gamma$ runs over the basis $\{T_i\}$
these sections form the rows of a fundamental matrix of flat sections,
and the $J$-function of $X$ is by definition
$$J^X:=\sum_{i=0}^s\left (\int_XF_{T_i}^X\right){\check T_i}
=\sum_{i=0}^s\langle\langle
\frac{T_i e^{(t_1T_1+\dots t_mT_m)/\hbar}}{\hbar-\psi},1
\rangle\rangle^X{\check T_i} .$$

We now consider these matrices corresponding to the standard basis
of $H^*({\bf G})$, and the basis ${\cal B}$ of $H^*({\bf P})$.
By Theorem 2.7, 
$$\left. \frac{p_V(F_{\sigma _{\nu} \Delta}^{\bf P} )}{\Delta}
\right|_{t_i = t + (r-1)\pi \sqrt{-1}} 
= \sigma _\nu e^{(t+(r-1)\pi \sqrt{-1})\sigma _1/\hbar} + o(e^t)$$ 
is a flat section for $\bf G$. However, we also have
$$F_{\sigma_\nu e^{(r-1)\pi\sqrt{-1}\sigma_1}}^{{\bf G}}
=\sigma_\nu e^{(t+(r-1)\pi \sqrt{-1})\sigma _1/\hbar}+o(e^t).$$  
By the uniqueness of flat sections with given intial term in the 
$e^t$-expansion we conclude that, for any two partitions $\mu,\nu$,

$$\left.\frac{(-1)^{{r \choose 2}}}{r!}\langle\langle
\frac{\sigma_{\nu}\Delta e^{\sum_it_ix_i/\hbar}}
{\hbar -\psi},\sigma_{\mu}\Delta\rangle\rangle^{\bf P}
\right|_{t_i=t+(r-1)\pi\sqrt{-1}}$$
\begin{equation}=\langle\langle
\frac{\sigma_{\nu}e^{(t+(r-1)\pi\sqrt{-1})\sigma_1/\hbar}}
{\hbar -\psi},\sigma_\mu\rangle\rangle^{{\bf G}}. \label{???}\end{equation}

On the other hand, the same argument as in the proof of 
Lemma 4.2 of \cite{JK}
shows that 
\begin{equation}
\int_{\bf P}F_{\sigma _{\nu}\Delta }^{{\bf P}}\cup\Delta
={\cal D}_{\Delta}\int_{\bf P}F_{\sigma _\nu\Delta }^{\bf P}
\label{??}
\end{equation} 
since, for any decomposition $\Delta=\prod (x_i-x_j) \prod (x_k-x_l)$,
there are no quantum correction terms
in the product $\prod (x_i-x_j) * \prod (x_k-x_l)$.

Next observe that multiplication by 
the class
$e^{\sigma_1(r-1)\pi\sqrt{-1}/\hbar}$ is an
invertible linear operator on $H^*({\bf G})$
(with inverse 
the operator of 
multiplication by
$e^{-\sigma_1(r-1)\pi\sqrt{-1}/\hbar}$), therefore
$$e^{\sigma_1(r-1)\pi\sqrt{-1}/\hbar}J^{{\bf G}}=\sum_{\nu}
\langle\langle\frac{\sigma_\nu e^{(t+(r-1)\pi\sqrt{-1})\sigma_1/\hbar}}
{\hbar -\psi},1
\rangle\rangle^{{\bf G}}\sigma_{\check{\nu}}.$$
It follows then from (\ref{???}) (applied to $\mu=$ the empty partition) that
$$e^{\sigma_1(r-1)\pi\sqrt{-1}/\hbar}J^{{\bf G}}=\left.
\sum_\nu \frac{(-1)^{{r \choose 2}}}{r!}\langle\langle
\frac{\sigma_\nu\Delta e^{\sum t_ix_i}}{\hbar -\psi},\Delta
\rangle\rangle^{\bf P} 
\sigma_{\check{\nu }}\right |_{t_i=t+(r-1)\pi\sqrt{-1}}$$
$$=\left. \sum_\nu \frac{(-1)^{{r \choose 2}}}{r!}\left(\int_{\bf P}
F_{\sigma_{\nu}\Delta}^{{\bf P}}\cup\Delta\right)\sigma_{\check{\nu}}
\right|_{t_i=t+(r-1)\pi\sqrt{-1}}.$$
From this and (\ref{??}) we get
$$\left. e^{\sigma_1(r-1)\pi\sqrt{-1}/\hbar}J^{{\bf G}}={\cal D}_{\Delta} 
\sum_\nu\frac{(-1)^{{r \choose 2}}}{r!}\left(\int_{\bf P}
F_{\sigma _\nu\Delta }^{\bf P}\right)\sigma _{\check{\nu }}  
\right|_{t_i=t+(r-1)\pi \sqrt{-1}}$$
$$=\frac{1}{\Delta}p_V \left(\left.{{\cal D}_{\Delta}}(J^{\bf P})  
\right|_{t_i=t +(r-1)\pi \sqrt{-1}}\right).$$
To finish the proof we simply note that
$$\frac{1}{\Delta} p_V \left(\left.{\cal D}_{\Delta}(J^{\bf P})  
\right|_{t_i=t +(r-1)\pi \sqrt{-1}}\right)=
\frac{1}{\Delta}{\cal D}_{\Delta}(J^{\bf P})   
|_{t_i=t + (r-1)\pi\sqrt{-1}}$$ 
since $\left.{\cal D}_{\Delta}(J^{\bf P})  
\right|_{t_i=t +(r-1)\pi \sqrt{-1}}$ is already antisymmetric.

\bigskip

\nt{\bf Remark:} Since multiplication by $e^{t\sigma _1}$ 
is an invertible operator on $H^*({\bf P})$, equation (\ref{???}) gives
$$\left. \frac{(-1)^{{r \choose 2}}}{r!}
\langle\langle\frac{\sigma_\nu\Delta }
{\hbar - \psi},\sigma_{\mu}\Delta\rangle\rangle^{\bf P}
\right|_{t_i = t + (r-1)\pi\sqrt{-1}}
=\langle\langle\frac{\sigma_\nu}{\hbar -\psi},\sigma_\mu
\rangle\rangle^{{\bf G}}, $$ for any two partitions $\mu, \nu$.
Just as the formula (\ref{qintegration}) in the case of 
three-point Gromov-Witten invariants,
this can be viewed as an analogue
of Martin's integration formula for two-point descendants.

\newpage

\nt {\bf \S 3. Applications.}
\bigskip

The first application we discuss is a proof (for Grassmannians)
of the ``$R$-conjecture'' of Givental. According to Givental's program
\cite{Giv3}, once the $R$-conjecture is established, the Virasoro conjecture
for Grassmannians follows. 
For this we will need 
the Hori-Vafa conjecture in its original formulation \cite{HV}, which we now
explain. 

There is an extension of Mirror Symmetry to Fano manifolds (see  
\cite{EHX1}, \cite{Giv1}, 
\cite{HV}). If $X$ is Fano, the general Physics formulation of it says
that the nonlinear sigma model for $X$ is equivalent to the 
Landau-Ginzburg theory
with a certain potential function $W$. 
For $X={\bf P}$
Hori and Vafa (\cite{HV}) argue on physical grounds
that the mirror symmetric Landau-Ginzburg theory is given by
$$W=\sum_{i=1}^rW_i,\;\; {\mathrm{with}}\;\; W_i=\sum_{j=1}^{n}y_{ij},$$
where the variables $y_{ij}$ satisfy the relations 
$\prod_{j=1}^ny_{ij}=e^{t_i}$, $i=1,\dots ,r$.

While the equivalence cannot be expressed in this generality in 
rigorous mathematical
terms, a precise statement about integral representations for
the coefficients of the $J$-function can be extracted from it. 
This statement has been 
formulated and proved by Givental, several years before \cite{HV} was written. 

\nt {\bf Theorem 3.1 (\cite{Giv1}):} For each $\hbar \neq 0$,
the coefficients of $J^{{\bf P}}$ -- with 
respect to every basis of $H^*({\bf P})$ -- span the same
${\bf C}$-subspace (of dimension $n^r={\mathrm rank}(H^*({\bf P})$) 
in the vector space 
${\bf C}[t_1,\dots ,t_r][[e^{t_1},\dots ,e^{t_r}]]$
as do the functions
$$I_K=\int_{\prod_{i=1}^n \Gamma_{i,k_i}}e^{W/\hbar}\frac{\wedge_{ij} dy_{ij}}
{\wedge_id(\prod_j y_{ij})}, \; \;
K=(k_1,\dots ,k_r), \;\; 1\leq k_i\leq n, $$
where, for each $i$, 
$\Gamma_{i,k_i}\subset Y_{t_i}:=\{ \prod_{j=1}^n y_{ij}=e^{t_i}\}$
is cycle representing a relative homology class in 
$H_{n-1}(Y_{t_i},{\mathrm{Re}}(W_i/\hbar)=-\infty)$.

\medskip

One gets a basis of this subspace by choosing for each $1\leq i\leq r$
the basis 
$\{\Gamma_{i,1},\dots , \Gamma_{i,n}\}$ of
$H_{n-1}(Y_{t_i},{\mathrm{Re}}(W_i/\hbar)=-\infty)$ consisting of the
descending Morse cycles corresponding to the $n$ non-degenerate
critical points of $W_i$. In what follows we will always consider
this choice of $\Gamma_{i,k_i}$'s, and denote $\Gamma_K:=
\prod_{i=1}^n \Gamma_{i,k_i}$.

\medskip

Furthermore, in the same paper Hori and Vafa conjecture that 
the Landau-Ginzburg
theory mirror-symmetric to the 
nonlinear sigma model on the Grassmannian ${\bf G}$ is obtained
by a kind of symmetrization
from the Landau-Ginzburg theory mirror to ${\bf P}$. 
As a consequence, they conjecture 

\medskip

\nt {\bf Theorem 3.2 (Hori-Vafa Conjecture - integral representation form):} 
For each $\hbar\neq 0$, the functions
$$\left.{\cal D}_{\Delta}I_K\right|_{t_i=t+(r-1)\pi\sqrt{-1}},\;\;\; \;\; 
K=(k_1,\dots ,k_r), \;\; 1\leq k_i\leq n$$
span the same 
${\bf C}$-subspace of ${\bf C}[t][[e^t]]$ as do the coefficients of
$J^{{\bf G}}$.

\medskip

{\bf Proof:} The $J$-function version proved in this paper gives
$$J^{{\bf G}}=\left. e^{-\sigma_1(r-1)\pi \sqrt{-1}/\hbar} 
\frac{{\cal D}_{\Delta}
J^{\bf P}}{\Delta}\right |_{t_i=t+(r-1)\pi\sqrt{-1}}.$$ 
Since multiplication by $e^{-\sigma_1(r-1)\pi \sqrt{-1}/\hbar}$ 
is an invertible linear transformation on $H^*({{\bf G}})$, it follows
that the overall factor $e^{-\sigma_1(r-1)\pi \sqrt{-1}/\hbar}$ 
can be discarded from $J^{{\bf G}}$ without changing the span of the 
coefficients. The statement follows 
now immediately from Theorem 3.1.

\medskip

\nt{\bf Remark:} Note that 
$\left.{\cal D}_{\Delta}I_K\right|_{t_i=t+(r-1)\pi\sqrt{-1}}$ vanishes 
whenever we have $k_i=k_j$ for some $i\neq j$. Therefore we can restrict 
in the statement of the theorem to those $K$'s with distinct entries.

\bigskip

Recall that ${{\bf G}}$ has a natural action of the torus 
$T=({\bf C}^*)^n$. The 
torus also acts diagonally on ${\bf P}$ and there are $T$-equivariant
versions of $J^{{\bf P}}$ and of Givental's theorem above.
Explicitly, it is shown in \cite{Giv1} that the formula
$$J^{{\bf P}^T}= 
e^{\frac{t_1x_1 + ... + t_rx_r}{\hbar}} \sum_{(d_1,...,d_r)} 
\frac{e^{d_1t_1 + ... + d_rt_r}}{\prod_{i=1}^r\prod_{l=1}^{d_i}
\prod _{j=1}^n (x_i-\lambda _j+ l\hbar )}$$
holds, and in \cite{Giv3} that, for fixed $\hbar\neq 0, \kl_i\neq 0$,
$i=1,\dots ,n$, 
the coefficients of $J^{{\bf P}^T}$ (in
any ${\bf C}[\lambda_1,\dots ,\lambda_r]$-basis of the equivariant
cohomology ring $H^*_T({\bf P})$) generate the same  
${\bf C}$-subspace of 
${\bf C}[t][[e^t]]$ as
do the integrals 
$$I_K^T=\int_{\prod \Gamma_{i,k_i}^T}e^{W^T/\hbar}\frac{\wedge_{ij} dy_{ij}}
{\wedge_id(\prod_j y_{ij})}, \; \;
K=(k_1,\dots ,k_r), \;\; 1\leq k_i\leq n. $$
The only difference between the integrals $I_K$ and $I_K^T$
is in the ``phase function'' $W^T$, which is now given by
$$W^T=\sum_{i=1}^rW_i^T,\;\; {\mathrm{with}}\;\; 
W_i^T=W_i+\sum_{j=1}^n\lambda_j\ln y_{ij}=
\sum_{j=1}^{n}y_{ij}+\lambda_j\ln y_{ij}.$$
Consequently, the cycles
$\Gamma_{i,k_i}^T\subset Y_{t_i}:=\{ \prod_{j=1}^n y_{ij}=e^{t_i}\}$
represent relative homology classes in 
$H_{n-1}(Y_{t_i},{\mathrm{Re}}(W_i^T/\hbar)=-\infty)$.  

The same argument as above, but using Theorem 1.5' gives

\medskip

\nt{\bf Theorem 3.2':} For fixed $\hbar\neq 0, \kl_i\neq 0$,
$i=1,\dots ,n$, the functions
$$\left.{\cal D}_{\Delta}I_K^T\right|_{t_i=t+(r-1)\pi\sqrt{-1}},\;\;\; \;\; 
K=(k_1,\dots ,k_r), \;\; 1\leq k_i\leq n$$
span the same 
${\bf C}$-subspace of 
${\bf C}[t][[e^t]]$ as do the 
coefficients of $J^{{\bf G}^T}$. 

\medskip

Theorems 3.2 and 3.2' are ``mirror theorems'' for the Grassmannian,
in the sense of Givental.

\bigskip

We now turn to Givental's $R$-Conjecture and first explain briefly the
statement, refering to \cite{Giv3}, \cite{JK} for details.

Let $X$ be a Fano
projective manifold with an action of a torus $T=({\bf C}^*)^n$.
As before, we denote by $\lambda_1,\dots ,\lambda_n$ the generators
of $H^*_T(BT)$. Assume that
there are $T$-equivariant nef divisor classes $\gamma_i$,
$i=1,\dots ,m$ forming a
basis of $H^2_T(X)$. As usual, introduce dual parameters $t_i$ and set
$q_i =e^{t_i}$. All structures discussed in the previous sections
of the paper have equivariant analogues:
the $T$-equivariant Gromov-Witten invariants of $X$ (\cite{Giv1}, \cite{Kim2}) 
can be used to construct a (small) equivariant quantum product $*^T$
on $H^*_T(X)\otimes {\bf C}[e^{t_1},\dots ,e^{t_m}]$,
and we have the associated flat connection $\nabla_{\hbar}$ 
defined by the same formula. 
Choose a ${\bf C}[\lambda_1,\dots ,\lambda_n]$-basis
$\{\phi _j\}$ of $H^*_T(X)$. A fundamental matrix of
flat sections is constructed using
two-point descendants, and its last column
gives the equivariant $S$-function.

Let $A_i$ be the matrix in the basis $\{\phi _j\}$ of the operator
on $QH^*_T(X)$ given by quantum multiplication
by $\gamma_i$.
Assume furthermore that the $T$-fixed points on $X$ are isolated and 
the $1$-dimensional orbits are isolated as well. Then the
equivariant cohomology ring $H^*_T(X)$, and therefore also
the equivariant small quantum
cohomology ring $QH^*_T(X)$ are semi-simple. In this case there
are {\it canonical coordinates} $u_j$ on $H^*_T(X)$,
unique up to addition of terms constant in the $t_i$'s,
such that the 
connection matrix of 1-forms matrix has 
$N:={\mathrm rank}(H^*(X))$ distinct
eigenvalues $du_j$ for general $\lambda$ and $q$.

Denote by $\Psi$ the eigenmatrix of $\sum_iA_i dt_i$, normalized
so that the row-vectors $\sum_j \Psi _{jk}\phi_j$ are unit
vectors with respect to the Poincar\'e metric. We then have
$$(\sum_i A_i dt_i) \Psi = \Psi dU, $$ where
$U={\mathrm diag} (u_1,...,u_N)$.

It is shown in \cite{Giv2}, \cite{Giv3} that there is an asymptotic fundamental
solution matrix of the connection $\nabla _{\hbar}$ of the form 
$\Psi Re^{U/\hbar}$, with 
$$R=1+R_1\hbar + R_2\hbar ^2 +\dots$$
and $R_i$ analytic in $q$ and $\lambda$ near any point
$(q,\lambda )$ at which $QH^*_T(X)$ is semi-simple
(i.e., $\hbar\frac{\partial}{\partial t_i} (\Psi R
e^{U/\hbar}) \sim A_i \Psi R e^{U/\hbar}$ asymptotically in $\hbar
^l e^{U/\hbar}$, $l\rightarrow \infty$). The solution is unique
if the following two conditions are imposed:

1) the orthogonal condition: $R^t(\hbar)R(-\hbar)=1$,

2) the classical limit condition: for generic values of $\kl_1,\dots ,\kl_n$
$$\lim_{q\rightarrow 0}R=e^ {{\mathrm diag} (b_1,...,b_N)},$$ where
$$ b_i(\hbar ) = \sum _{k=1}^{\infty} N_{2k-1}^{(i)}
\frac{B_{2k}}{2k} \frac{\hbar ^{2k-1}}{2k-1}.$$ 
Here $B_{2k}$ are
Bernoulli numbers defined by 
$$x/(1-e^{-x}) = 1+x/2+ \sum_{k=1}^{\infty} B_{2k}\frac{x^{2k}}{(2k)!}.$$ 
Finally, the $N_l^{(i)}$ are defined
as follows: There is a natural 1-1 correspondence between the $T$-fixed
points $\{v_i\; |\; i=1,\dots ,N\}$ in $X$ and  
idempotents $\frac{\partial}{\partial u_i}$ in 
$H^*_T(X)\otimes_{{\bf C}[\kl_1,\dots ,\kl_n]}{\bf C}
(\lambda _1,...,\lambda _n)$. Let $w_{ij}$ denote the weights of 
the $T$ action on the tangent spaces $T_{v_i}X$. With this notation, 
$$N_l^{(i)} = \sum _j\frac{1}{w_{ij}^l}.$$ 

\medskip

\nt{\bf Remark:} In \cite{JK} an additional requirement (the ``equivariant
homogeneity condition'') is listed as necessary to get the uniqueness of $R$.
In fact, this condition is always satisfied once 1) and 2) above are imposed,
as can be seen from the inductive construction of $R$ in \cite{Giv2}. 

\bigskip

We are ready to state Givental's $R$-conjecture (in fact, the
``$R$-conjecture restricted to $H^2$'', in the terminology of \cite{JK}).

\medskip

\nt{\bf R-Conjecture (\cite{Giv3}):} The non-equivariant limit 
$\lambda_1\rightarrow 0,..., \lambda_n\rightarrow 0$ of the
matrix $R$ exists.

\bigskip

Givental \cite{Giv3} proves the conjecture in the case $X$ 
is a projective space,
and in fact the same proof works in the case of a product of projective
spaces. In particular, the conjecture holds for our ${\bf P}$.
The main ingredient of the proof is the existence of a fundamental solution
matrix of the flat connection whose entries are the integrals from 
Theorem 3.1 and their derivatives.

The validity of the $R$-conjecture for ${\bf G}$ would follow immediately 
from its validity for ${\bf P}$ if
we would have at our disposal the equivariant version of Theorem 2.7.
In turn, this would follow from an equivariant version of the ``quantum
integration formula'' (\ref{qintegration}). However, neither of the two
proofs we gave to this formula are readily generalized. Indeed, there
is no analogue of the rim-hook algorithm for equivariant quantum
cohomolgy, while to establish the equivariant Vafa-Intriligator formula
one needs an equivariant version of the quantum Giambelli formula
of \cite{Ber2}, which is not known. In short, although we believe that the 
equivariant quantum integration formula is true, lacking this, 
we give a direct proof to the $R$-conjecture for ${\bf G}$
which uses Theorem 3.2' and parallels
the proofs in \cite{Giv3} and \cite{JK}.  

\medskip

\nt{\bf Theorem 3.3:} The R-conjecture holds for ${\bf G}=G(r,n)$.

{\bf Proof:}
We start by noting that it follows from the localization theorem
that the powers 
$$\{\sigma_1^i\; |\;  i=0,1,...,{\mathrm rank}(H^*({\bf G}))-1\}$$
form a basis over ${\bf C}(\kl_1,...,\kl_n)$
for the localized equivariant cohomology ring 
$$H^*_T({\bf G})\ot_{{\bf C}[\kl_1,...,\kl_n]}{\bf C}(\kl_1,...,\kl_n)$$  
(cf. the proof of Lemma 1.3).

By Theorem 3.2' and Lemma 4.2 of \cite{JK}, there exist
differential operators $D_i,\; i=1,\dots ,N$, in 
${\bf C}(\lambda_1,\dots ,\lambda_n)[e^t,\hbar ,\hbar\frac{d}{dt}]$ such that
$$s_K:=\sum_i \sigma_1^i D_i 
({\cal D}_{\Delta} I_{K}^T|_{t_j=t+(r-1)\pi\sqrt{-1}}).$$
form a fundamental solution to the flat
connection associated to the equivariant small quantum cohomology of ${\bf G}$.
By clearing denominators, we may assume that
the differential operators are in fact in 
${\bf C}[\lambda_1,\dots ,\lambda_n,e^t,\hbar ,\hbar\frac{d}{dt}]$. 
Now write $\sigma _1 ^i$ as a linear
combination of $\phi_j$ over ${\bf C} [\kl_1,...,\kl_n]$. It follows that
$$s_K=\sum_i s_{K,i}\phi_i=
\sum_i\phi_i\left.\left(\int_{\Gamma _K^T} e^{W^T/\hbar}
\varphi_{i,t,\hbar, \{\kl_j\}}
\frac{\wedge_{ij} dy_{ij}}
{\wedge_id(\prod_j y_{ij})}\right)\right |_{t_j=t+(r-1)\pi\sqrt{-1}}$$ 
form a fundamental solution, where the functions 
$\varphi_{i,t,\hbar, \{\kl_j\}}$ 
depend analytically on $q=e^t$, $\hbar$, and the $\lambda_j$'s.
The stationary phase approximation of the integrals
$s_{K,i}$ near the non-degenerate critical point corresponding to the cycle
$\Gamma_K^T$ gives an asymptotic expansion of the form
$$(s_{K,i}) \sim \Psi ' R e^{U/\hbar },$$
for some matrices $\Psi '$ and $R=1+R_1\hbar + R_2\hbar ^2 +\dots$

Since the integrand depends analytically on the $\kl_j$'s 
(and the usual small quantum cohomology ring $QH^*({\bf G})$
is semi-simple \cite{Abr}), it follows that
for the matrix $R$ constructed this way  
the limit $\lim_{\kl_j\rightarrow 0}R$ exists; 
moreover, since the $s_{K,i}$'s are the entries of 
a fundamental solution to the flat
connection, one can show that the eigenmatrices $\Psi '$ and 
$\Psi$ are related by $\Psi ' = C\Psi$, where $C$ is a
matrix constant with respect to $t$ (see e.g., the proof of Theorem 4.4
in \cite{JK} for details).
Therefore $\Psi R e^{U/\hbar}$ is a fundamental solution matrix.
It is now routine to check that
the orthogonality and classical limit conditions hold for $R$ 
(cf. \cite{Giv3}, \cite{JK}).
Hence the matrix $R$ constructed in this manner is indeed the matrix
in the $R$-conjecture.

\bigskip

According to \cite{Giv3}, if the (non-equivariant) small 
quantum cohomology ring $QH^*(X)$
is also semi-simple, then the Virasoro Conjecture of Eguchi-Hori-Xiong
\cite{EHX2}, and
S. Katz follows from the $R$-conjecture. Therefore, Theorem 3.3 has the
following 

\medskip

\nt{\bf Corollary 3.4:} The Virasoro conjecture is true for ${\bf G}$.

\bigskip

As our second application, we present a proof for $G(2,n)$ of a conjecture
made in \cite{BCKS1}. In that paper,
another ``mirror construction'' for Grassmannians (indeed, for all 
partial flag manifolds, see \cite{BCKS2})
was proposed. The construction is
via toric degenerations, and leads to
integral representations for the coefficients of the $J$-function
that are completely different than the ones in Theorem 3.2 above.
As a consequence, an explicit formula for the coefficient of 
the cohomology class $1$
in $J^{G(r,n)}$ was conjectured. 
We were able at the time to verify it by brute force computation only
in a few cases ($G(2,n)$, $n\leq 7$, and $G(3,6)$) 
and then used it, together with the ``quantum Lefschetz hyperplane theorem'' 
to calculate
the genus $0$ Gromov-Witten invariants of some complete intersection 
Calabi-Yau $3$-folds.

Since in this paper we give a formula for the full $J$-function,
the proof of the conjecture reduces to checking the following 
(highly nontrivial) combinatorial identity:

\medskip

\nt{\bf Conjecture :} The constant term with respect to the $x_i$'s from 
$J^{G(r,n)}$ in Theorem 1.5 equals $A_{G(r,n)}(e^t)$ in \cite{BCKS1}, 
Conjecture 5.2.3.

\medskip

We now make this more explicit in the case $r=2$. First, the formula
in \cite{BCKS1} is

$$A_{G(2,n)}(e^t)=\sum_{d\geq 0}\frac{e^{dt}}{(d!)^n}
\sum_{j_{n-3}\geq ...\geq j_1}{d\choose j_{n-3}}^2
{d\choose j_{n-4}}...{d\choose j_1}
{j_{n-3}\choose j_{n-4}}...{j_2\choose j_1},$$
which may be rewritten 
$$
\sum_{d\geq 0}\frac{e^{dt}}{(d!)^n}
\sum_{j_{n-3}\ge j_{n-4}\ge \cdots \ge j_1\ge j_0=0} 
\frac{d!^{n-2}}{(d-j_{n-3})!\prod_{i=2}^{n-2} (d-j_{i-1})! 
(j_{i-1}-j_{i-2})!j_{i-1}!}.$$

On the other hand,
extracting the constant term with respect to the $x_i$'s from $J^{G(2,n)}$ 
(and setting $\hbar =1$) in 
Theorem 1.5 we get
$$\sum_{d\geq 0}\frac{e^{dt}}{(d!)^n}
\frac{(-1)^{d}}{2} \sum_{m=0}^d {d\choose m}^{n}
\biggl( n(d-2m)(\gamma(m)-\gamma(d-m))+2\biggl),
$$ 
where for a positive integer $m$ 
$$
\gamma(m)=\sum_{j=1}^m \frac{1}{j}.
$$
is the partial sum of the harmonic series, and $\gamma(0)=0$. 

The conjecture for $G(2,n)$ follows
then from

\nt{\bf Proposition 3.5:} If $d$ are $n$ are nonnegative 
integers, $n\geq 3$, then
$$
\sum_{j_{n-3}\ge j_{n-4}\ge \cdots \ge j_1\ge j_0=0} 
\frac{d!^{n-2}}{(d-j_{n-3})!\prod_{i=2}^{n-2} (d-j_{i-1})! 
(j_{i-1}-j_{i-2})!j_{i-1}!}$$
$$=\frac{(-1)^{d}}{2} \sum_{m=0}^d {d\choose m}^{n}
\biggl( n(d-2m)(\gamma(m)-\gamma(d-m))+2\biggl).
$$ 

\medskip

{\it Proof (Dennis Stanton):} A formula expressing a $k$-fold sum as 
a sum over a single summation
index is given by the 
iterate of Bailey's lemma. Explicitly, if in \cite{And}, p. 30, Theorem 3.4
we change notation by setting $n=m$, $N=d$, $k=n-2$, and $n_i=j_{i-1}$, then we choose
 
$$
\beta_m=
\left\{ \begin{array}{ll} 
1 & {\mathrm if}\quad  m=0 \\
0 & {\mathrm if}\quad  m>0,
\end{array}\right.
\quad
\alpha_m=(-1)^m q^{m\choose 2}\frac{(a;q)_m}{(q;q)_m}
\frac{1-aq^{2m}}{1-a},
$$
all $c_i=q^{-d}$, and let all $b_i\to \infty$, the identity in  
{\it loc. cit.} becomes
$$
\frac{(aq;q)_d}{ (aq^{1+d};q)_d}
\sum_{j_{n-3}\ge j_{n-4}\ge \cdots \ge j_1\ge j_0=0}
\frac{(q^{-d};q)_{j_{n-3}}\prod_{i=2}^{n-2} (q^{-d};q)_{j_{i-1}}}
{\prod_{i=2}^{n-2} (aq^{1+d};q)_{j_{i-1}}(q;q)_{{j_{i-1}}-j_{i-2}}}$$
$$\times
(q^{-2d}/a)^{j_{n-3}} (-1)^{j_{n-3}} q^{-{j_{n-3}\choose 2}}
\prod_{i=2}^{n-2} (-1)^{j_{i-1}}  q^{dj_{i-1}}a^{j_{i-1}}q^{j_{i-1}+1\choose 2} =$$
$$\sum_{m=0}^d \frac{(q^{-d};q)_m^{n-1}}{(aq^{1+d};q)_m^{n-1}}
(-1)^{(n-1)m}q^{(n-2){m\choose 2}}a^{(n-2)m}q^{(n-2+d)m+dm(n-2)} q^{-{m\choose 2}}
\alpha_m,
$$
where for a nonnegative integer $m$ 
$$(a;q)_m=\prod_{l=0}^{\infty}\frac{(1-aq^l)}{(1-aq^{l+m})}.$$

The $q\to 1$ limit of this equation when $a=q^{-d}$ is
$$
\sum_{j_{n-3}\ge j_{n-4}\ge \cdots \ge j_1\ge j_0=0} 
\frac{d!^{n-2}}{(d-j_{n-3})!\prod_{i=2}^{n-2} (d-j_{i-1})! (j_{i-1}-j_{i-2})!j_{i-1}!}=$$
$$\lim_{a\to -d}  \frac{d!}{(a+1)_d}
\sum_{m=0}^d {d\choose m}^{n-1}\frac{m!^{n-1}}{(a+1+d)_m^{n-1}}
\frac{(a)_m}{m!}\frac{a+2m}{a}(-1)^m,$$
where $(a)_m$ is the generalized factorial 
$$(a)_m=a(a+1)...(a+m-1).$$

Notice that the left-hand side of the last identity coincides with the
left hand-side of the identity in the Proposition, while the right-hand side
is equal to
$$(*)\;\;\;\;\;\;\;\;\;\;\;\;\;\;\;\;\;\;\;\;\;\;\;\;\;\;\;\;\;\;\;\;\;\;\;\;\;\;\;\; 
d(-1)^{d-1} F'(-d),\;\;\;\;\;\;\;\;\;\;\; \quad\quad\quad
$$
where
$$
F(a)=\sum_{m=0}^d{d\choose m}^{n-1}\frac{m!^{n-1}}{(a+1+d)_m^{n-1}}
\frac{(a)_m}{m!}\frac{a+2m}{a}(-1)^m.
$$

However,
$$(**)\;\;\;
F'(-d)=\sum_{m=0}^d{d\choose m}^{n-1}
\biggl( -(n-1)\gamma(m){d\choose m} \frac{-d+2m}{-d}+$$
$${d\choose m}(\gamma(d-m)-\gamma(d))\frac{-d+2m}{-d}
+{d\choose m}\frac{-2m}{d^2}\biggr)=$$
$$
\quad\frac{-1}{d}\sum_{m=0}^d{d\choose m}^{n}
\biggl( \bigl((n-1)\gamma(m)+\gamma(d)-\gamma(d-m)\bigr)(d-2m)+2m/d
\biggr)
$$

Replacing $m$ by $d-m$ in this sum and adding to $(**)$ yields
$$
F'(-d)=\frac{-1}{2d}\sum_{m=0}^d{d\choose m}^{n}
\biggl( n(\gamma(m)-\gamma(d-m))(d-2m)+2\biggr).
$$

If we substitute now this expression for $F'(-d)$ into $(*)$, the
right-hand side of the identity in the Proposition is obtained.
This finishes the proof.

\medskip
\noindent Department of Mathematics, University of Utah, Salt Lake
City, UT 84112, {\it bertram@math.utah.edu}
\medskip

\noindent School of Mathematics, University of Minnesota,
Minneapolis MN, 55455,\\ {\it ciocan@math.umn.edu}
\medskip 

\noindent
Department of Mathematics, Pohang University of Science and Technology,\\
Pohang, 790-784, Republic of Korea, {\it bumsig@postech.edu}

\end{document}